\def\AA         {{\bf A}}
\def\ZZ         {{\bf Z}}
\def\RR         {{\bf R}}
\def\CC         {{\bf C}}
\def\QQ         {{\bf Q}}
\def\PP         {{\bf P}}

\def\ee		{{\rm e}}
\def\no		{:}
\def\reg	{{\rm reg.}}
\def\vac	{\hspace{-4pt}>}

\def\Fock	{{\rm Fock}}
\def\Bfg	{{\cal BRST}_{f,g}}
\def\MSV	{chiral de Rham complex}
\def\msv	{{\cal MSV}}
\def\VOA	{{\cal VA}}
\def\dim	{{\rm dim}}
\def\BRST	{{\cal BRST}}

\documentstyle[twoside,12pt]{article}
\newtheorem{prop}{Proposition}[section]
\newtheorem{dfn}[prop]{Definition}
\newtheorem{theo}[prop]{Theorem}
\newtheorem{conj}[prop]{Conjecture}
\newtheorem{rem}[prop]{Remark}
\newtheorem{coro}[prop]{Corollary}
\newtheorem{lem}[prop]{Lemma}

\title{Vertex Algebras and Mirror Symmetry}

\author{
Lev A. Borisov
\\
\small Department of Mathematics,  Columbia University \\
\small 2990 Broadway, Mailcode 4432, New York, NY 10027, USA \\
\small e-mail: lborisov@math.columbia.edu}

\begin{document}

\date{}

\maketitle

\begin{abstract}
{
Mirror Symmetry for Calabi-Yau hypersurfaces in toric
varieties is by now well established. However, previous
approaches to it did not uncover the underlying reason 
for mirror varieties to be mirror. We are able to calculate 
explicitly vertex algebras that correspond to holomorphic 
parts of A and B models of Calabi-Yau hypersurfaces
and complete intersections in toric varieties. We establish
the relation between these vertex algebras for mirror Calabi-Yau
manifolds. This should eventually allow us to rewrite the whole story of
toric Mirror Symmetry in the language of sheaves of vertex algebras. Our
approach is purely algebraic and involves simple techniques 
from toric geometry and homological algebra, as well as 
some basic results of the theory of vertex algebras. Ideas of this
paper may also be useful in other problems related to maps
from curves to algebraic varieties.

This paper could also be of interest to physicists, because it contains
explicit description of holomorphic parts of A and B models of Calabi-Yau
hypersurfaces and complete intersections in terms of free bosons and
fermions.
}
\end{abstract}

\section{Introduction}
First example of Mirror Symmetry was discovered by physicists
in \cite{candelasetal}. It relates A model on one Calabi-Yau
variety with B model on another one. Unfortunately, the definition
of A and B models was given by physicists in terms of integrals over
the set of all maps from Riemann surfaces to a given Calabi-Yau
variety, see \cite{witten}. While physicists have developed good
intuitive understanding of the behavior of these integrals, they are
ill-defined mathematically. Nevertheless, physicists came up with
predictions of numbers of rational curves of given degree in
Calabi-Yau manifolds, a quintic threefold being most prominent example.

Kontsevich has introduced spaces of stable maps (see \cite{kontsevich})
which allowed him  to define mathematically virtual numbers of rational
curves on a quintic. Givental proved in \cite{givental}
that these virtual numbers agree with physical predictions. Because
of its hard calculations, Givental's paper is a source of controversy, see
\cite{yauetal} by Lian, Liu and Yau. In some sense, however, Givental's
approach does not clarify the origins of mirror symmetry. His proof is a
beautiful and tricky calculation which has little to do with mirror
involution.

The goal of this paper is to present a completely different approach
to toric Mirror Symmetry which should eventually lead to conceptual
understanding of mirror involution in purely mathematical terms.
To do this, one has to employ the theory of vertex algebras, which is
a very well developed purely algebraic theory.
Malikov, Schechtman and Vaintrob have recently suggested an algebraic
approach to A models (see  \cite{MSV}) that involves \MSV\
which is a certain sheaf of vertex  algebras. In my personal
opinion, their paper is one of the most important mathematics papers on mirror
symmetry written to date, even though it does not deal directly with
mirror symmetry.

In this paper we are able to calculate A and B model vertex algebras
for mirror families of Calabi-Yau complete intersections. We define 
(quasi-)loop-coherent sheaves over any algebraic variety $X$, and
we show that if their sections on any affine open subset have 
vertex algebra structure, then the cohomology of the sheaf has
this structure as well. It is hoped that the techniques of this paper
will prove to be more important than the paper itself, after all,
they should allow mathematicians to do rigorously what physicists
have been doing half-rigorously for quite a while and with a lot
of success. 

As far as applications to conformal field theory are concerned, 
this paper suggests a way of defining A and B models for
varieties with Gorenstein toroidal singularities that does not
use any resolutions of such singularities. Also, we describe
explicitly vertex algebras that correspond to points in fully enlarged
K\"{a}hler cone (see \cite{greene}) and do not come from Calabi-Yau
varieties.

The paper is organized as follows. Section 2 is devoted to 
(quasi-)loop-coherent sheaves, which is a generalization
of the notion of (quasi-)coherent sheaves. It serves as a useful
framework for the whole paper, and perhaps may have other applications.
Sections 3 and 4 contain mostly background material. The only
apparently new result there is Proposition \ref{vertcoh}, which 
had actually been suggested in \cite{MSV}. Section 5 contains
a calculation of \MSV\ of a hypersurface in a smooth variety 
in terms of \MSV\ of the corresponding line bundle. Section 6
contains a calculation which is in a sense a mirror image of the
calculation of Section 5. Sections 7 and 8 combine results of previous
sections to describe A and B models of Calabi-Yau hypersurfaces
in smooth nef-Fano toric varieties. Sections 9 and 10 attempt
to extend these results to singular varieties and complete intersections.
While good progress is made there, a few details need to be 
further clarified. Section 11 is largely a speculation. We state 
there some open questions related to our construction, as well as some
possible applications of the results and techniques of this paper.

We use the book of Kac \cite{Kac} as the standard reference for
vertex algebras.

{\bf Acknowledgments.} This project began in 1995
while I was a Sloan predoctoral fellow at the University of
Michigan. I thank Martin Halpern and Christoph Schweigert
who helped me learn the basics of conformal field theory. 
Konstantin Styrkas has answered a couple of my questions
regarding vertex algebras, which has greatly improved my 
understanding of the subject. 

\section{Loop-coherent sheaves}
The goal of this section is to develop foundations of
the theory of (quasi-)loop-coherent sheaves over algebraic varieties. 
These are rather peculiar objects that nevertheless behave very much like
usual coherent and quasi-coherent sheaves. For simplicity, we 
only concern ourselves with algebras over complex numbers.
This is mostly just a formalization of the localization calculation of
\cite{MSV} but it provides us with a nice framework for our discussion.
The idea is somehow to work with sheaves over the loop space of an
algebraic variety without worrying much about infinities. Only future
can tell if this is a truly useful concept or just an annoying 
technicality.

\begin{dfn}
{\rm Let $R$ be a commutative algebra over $\CC$ with a unit. {\it
$R$-loop-module}\ is a vector space $V$ over $\CC$ together with the
following
set of data. First of all, $V$ is graded
$$V=\oplus_{l\geq 0}V_l.$$
We assume
that values of $l$ are integer, although it takes little
effort to modify our definitions to allow any real $l$. 
We denote by $L[0]$ the grading operator, that is $L[0]v=kv$
for all $v\in V_k$. 
In addition, for every element $r\in R$ and every integer 
$l$ there given a linear operator  $r[l]:V\to V$ such that the following
conditions hold\\
(a) $1[k]=\delta_k^0$ ;\\
(b) all $r[l]$ commute with each other;\\
(c) $r[l]V_k\subseteq V_{k-l}$ ;\\
(d) for every two ring elements $r_1$, $r_2$ there holds 
$$(\sum_k r_1[k]z^{-k})(\sum_l r_2[l]z^{-l})
=\sum_k (r_1r_2)[k]z^{-k}.$$
This equation makes sense because at any given power
of $z$ while applied to any given element, only the finite number of
terms on the left hand side are non-zero. This follows from (b), (c),
and $L[0]\geq 0$.
}
\end{dfn}

\begin{rem}
{\rm 
$R$-loop-modules are usually not $R$-modules. Really, one has 
$$(r_1r_2)[0]=r_1[0]r_2[0]+\sum_{k\neq 0}r_1[k]r_2[-k]$$
as opposed to just $(r_1r_2)[0]=r_1[0]r_2[0]$. However, 
the extra term is locally nilpotent, that is a sufficient 
power of it annihilates $V_l$ for any given $l$. This is
what makes it possible to localize loop-modules analogously
to localization of usual modules.
}
\label{nilp}
\end{rem}

The following proposition will also serve as a definition.
\begin{prop}
{\rm Let $S$ be a multiplicative system in $R$. Given a loop-module
$V$ over $R$ denote by $V_S$ its localization by the multiplicative system
$S_{loop}$ generated by $s[0]$ for all $s\in S$. We claim that $V_S$
has a natural structure of $R_S$-loop-module. Moreover, $\rho:V\to V_S$
is a universal morphism in the sense that for any map to
$R_S$-loop-module $\rho_1:V\to V_1$ that is compatible with $R\to R_S$
there exists a unique $R_S$-loop-module map $\rho_2:V_S\to V_1$
such that $\rho_1=\rho_2\circ \rho$.
}
\label{loc}
\end{prop}

{\em Proof.} First of all, let us provide $V_S$ with $R_S$-loop-module
structure. Grading is clearly unaffected by localization. For
every $v\in V_S$ we have $v=w/\prod s_i[0]$ and every $a\in R_S$ looks
like $a=b/\prod s_j$. We define
$$a(z)v={1\over {\prod s_i[0]}} a(z) w$$
$$a(z)w=b(z)(\prod{1\over {s_j}}(z)) w$$
$${1\over {s_j}}(z) w = {1\over {s_j[0]+\sum_{k\neq 0}s_j[k]z^{-k}}}w
=\sum_{l=0}^\infty (-1)^ls_j[0]^{-l-1}(\sum_{k\neq 0} s_j[k]z^{-k})^lw$$
which gives an element in $V_S[z,z^{-1}]$. Even though it seems that
infinite sums appear when one
applies $(1/s_j(z))$ several times, for any given $w$ most $s_j[k]$ with
positive $k$ could be safely ignored, since they annihilate $w$ and 
could always be pushed through other $s_*[*]$ by commutativity.

One, of course, needs to check that the above definition is
self-consistent.
It is clear that if you change $v$ to $s[0]v/s[0]$, the result 
stays the same. Changing $a$ to $sa/s$ requires a 
certain calculation, but does not present any major difficulties either.

For every map $V\to V_1$ where $V_1$ is $R_S$-loop-module notice that
$s[0]$ is invertible on $V_1$ for all $s\in S$. Really,
$$s[0](s^{-1})[0]={\bf 1}+{\rm locally~nilpotent}$$
so $s[0]$ is invertible by the same trick. Therefore, the map 
$V\to V_1$ can be naturally pushed through $V_S$. 
\hfill$\Box$

\begin{rem}
{\rm 
Because of Remark \ref{nilp}, it is enough to localize by $s[0]$
for those $s$ that generate $S$. 
}
\label{genloc}
\end{rem}

Proposition \ref{loc} allows us to define a {\it quasi-loop-coherent
sheaf}\ on any complex algebraic variety $X$ as follows. 

\begin{dfn}
{\rm
A sheaf ${\cal V}$ of vector spaces over $\CC$ is called {\it
quasi-loop-coherent}\ if for every affine subset ${\rm Spec}(R)\subset
X$  sections $\Gamma({\rm Spec}(R),{\cal V})$ form an
$R$-loop-module and restriction maps are precisely localization maps of
Proposition \ref{loc}.
}
\label{qlc}
\end{dfn}

Many results about quasi-loop-coherent sheaves could be deduced 
from the standard results about quasi-coherent sheaves due to the 
following proposition.

\begin{prop} 
{\rm
For every $R$-loop-module $V$ consider the following filtration
$$F^lV=\sum_{i,s_1,...,s_i,k_1,...,k_i}\prod_is_i[k_i]V_{\leq l}.$$
We  have $F^0V\subseteq F^1V\subseteq ...$, and 
$F^{l+1}V/F^lV$ has a natural structure of $R$-module. Moreover,
this filtration commutes with localizations.
}
\label{filtr}
\end{prop}

{\em Proof.} Locally nilpotent operators $(s_1s_2)[0]-s_1[0]s_2[0]$
push $F^{l+1}V$ to $F^lV$, which provides the quotient with the
structure of $R$-module. Filtration commutes with localization,
because $s[0]$ commute with $L[0]$. \hfill$\Box$

As a result of this proposition every quasi-loop-coherent sheaf
${\cal V}$ is filtered by other quasi-loop-coherent sheaves $F^l{\cal V}$ 
and all quotients are quasi-coherent.
It is also worth mentioning that the above filtration is finite on every
${\cal V}_k$ which prompts the following definition.

\begin{dfn}
{\rm
A quasi-loop-coherent sheaf is called loop-coherent, or {\it loco}\
if quasi-coherent sheaves $F^{l+1}{\cal V}\cap{\cal V}_k/
F^l{\cal V}\cap{\cal V}_k$ are coherent for all $k$ and $l$.
}
\label{locosheaves}
\end{dfn}

From now on we also use abbreviation {\it quasi-loco}\ in
place of quasi-loop-coherent.

\begin{rem}
{\rm 
A zero component ${\cal V}_0$ 
of a (quasi-)loco sheaf is (quasi-)coherent.
}
\label{zerocoh}
\end{rem}

\begin{prop}
{\rm 
For any affine variety $X$ and quasi-loco sheaf ${\cal V}$ on
it cohomology spaces $H^i(X,{\cal V})$ are zero for $i\geq 1$.
For any projective variety $X$ all cohomology groups of a loco sheaf are
finite dimensional for each eigen-value of $L[0]$. 
}
\label{acyclic}
\end{prop}

{\em Proof.} For both statements, one considers a specific eigen-value $k$
of $L[0]$ and then applies an induction on $l$ in $F^l{\cal V}\cap {\cal
V}_k$. \hfill$\Box$

\begin{rem}
{\rm
There is a one-to-one functorial correspondence between $R$-loop-modules
and quasi-loco sheaves over ${\rm Spec} R$.
}
\label{msh}
\end{rem}

\begin{rem}
{\rm
One can modify the definition of quasi-loco sheaves to allow
negative eigen-values of $L[0]$, as long as there is some bound
$L[0]\geq -N$ on them.
}
\end{rem}

\section{Sheaves of vertex algebras}
We follow \cite{Kac} in our definition of a vertex algebra.
We only consider vertex algebras over $\CC$.  

\begin{dfn} 
{\rm
(\cite{Kac})
A vertex algebra $V$ is first of all a super vector space over ${\bf C}$,
that is  $V=V_0\oplus V_1$ where elements of $V_0$ are called bosonic or
even
and elements from $V_1$ are called fermionic or odd. In addition,
there given a bosonic vector $|0\vac$ called {\it vacuum vector}\ .
The last part of the data that defines a vertex algebra is the
so-called state-field correspondence 
which is a parity preserving linear map from $V$ to ${\rm End}
V[[z,z^{-1}]]$
$$a\,\line(0,1){5}\hspace{-3.6pt}\to Y(a,z)=\sum_{n\in Z}a_{(n)}z^{-n-1}$$
that satisfies the following axioms

\noindent$\bullet${\bf translation covariance:}
$\{T,Y(a,z)\}_-=\partial_z Y(a,z)$ where $\{,\}_-$ denotes
the usual commutator and $T$ is defined by $T(a)=a_{(-2)}|0\vac$;

\noindent$\bullet${\bf vacuum:} $Y(|0\vac,z)={\bf 1}_V,~
Y(a,z)|0\vac_{z=0}=a$;

\noindent$\bullet${\bf locality:} $(z-w)^N\{Y(a,z),Y(b,w)\}_{\mp}=0$
for all sufficiently big $N$, where $\mp$ is $+$ if and only if 
both $a$ and $b$ are fermionic. The equality is understood as an identity
of formal power series in $z$ and $w$.
}
\end{dfn}

We will often write $a(z)$ instead of $Y(a,z)$. Linear operators 
$a_{(k)}$ will be referred to as {\it modes}\ of $a$.
Every two fields $a(z)$, $b(w)$ of a vertex algebra have {\it operator
product
expansion}
$$
a(z)b(w)=\sum_{i=1}^N {c^i(w)\over {(z-w)^i}}+
\no a(z)b(w)\no$$
where the meaning of the symbols $1\over {(z-w)^i}$ 
and $\no\no$ in the above formulas is
clarified in Chapter 2 of \cite{Kac}. 
We only remark here that operator product expansion contains information
about all super-commutators of the modes of $a$ and $b$, and the sum on
the right hand side is finite due to the locality axiom. The sum is called
the singular part and the $\no\no$ term is called the regular part.
When there is no singular part, the OPE is called {\it non-singular}\ and
it
means that all modes of the two fields in question super-commute.

In this paper we will only use graded vertex algebras.
\begin{dfn} 
{\rm
(\cite{Kac})
A vertex algebra $V$ is called {\it graded} if there is given an even
diagonalizable operator $H$ on $V$ such that 
$$\{H,Y(a,z)\}_-=z\partial_zY(a,z)+Y(Ha,z).$$
}
\end{dfn}

When $Ha=ha$ we rewrite $Y(a,z)$ in the form
$$Y(a,z)=\sum_{n\in{-h+Z}}a[n]z^{-n-h}$$
and again call $a[n]$ modes of $a$. Number $h$ is usually called the
{\it conformal weight}\ of $a$ or $a(z)$. We use
brackets to denote modes,
which differs from the notations of \cite{Kac} and \cite{MSV}.
Modes can hardly be confused with commutators, but we use
$\{,\}_\pm$ notation for the latter just in case.

\begin{dfn} 
{\rm
(\cite{Kac})
A graded vertex algebra $V$ is called {\it conformal}\ if there 
chosen an even vector $v$ such that the corresponding field
$Y(a,z)=L(z)=\sum_n L[n]z^{-n-2}$ satisfies 
$$L(z)L(w)={(c/2)\over(z-w)^4}+{2L(w)\over(z-w)^2}+{\partial L(w)
\over z-w}+{\rm reg.}$$
We also require $L[-1]=T$ and $L[0]=H$. Number $c$ here is
called
{\it central charge} or rank of the algebra.
}
\end{dfn}

We now combine the theories of quasi-loco sheaves and vertex 
algebras to define sheaves of graded vertex algebras over 
an algebraic variety $X$. Abusing the notations,
we denote the grading operator $H$ by $L[0]$ even if
algebra $V$ is not conformal.

\begin{dfn}
{\rm
Let $R$ be a commutative algebra over $\CC$.
A graded vertex algebra $V$ is called vertex $R$-algebra
if $R$ is mapped to $L[0]=0$ component of $V$ 
so that images of all $r\in R$ are bosonic, all modes
$r[n]$ commute with each other and
$$Y(r_1,z)Y(r_2,z)=Y(r_1r_2,z).$$
In addition, we assume that $L[0]$ has only non-negative 
eigen-values.
}
\end{dfn}

Definition of graded algebra implies that 
$r[n]$ decreases eigen-values of $L[0]$ by $n$. Thus
any vertex $R$-algebra has a structure of $R$-loop-module.

\begin{prop}
{\rm 
Let $S$ be a multiplicative system in $R$
and let $V$ be a vertex $R$-algebra. Then 
the localization $V_S$ defined in Proposition 
\ref{loc} has a natural structure of vertex $R_S$-algebra.
}
\label{vertloc}
\end{prop}

{\em Proof.} For any element $a\in V$ and any set of elements
$s_i\in S$ we need to define the field $Y(a/\prod s_i[0],z)$
on $V_S$. First, we do it for any element of the 
form $|0\vac\hspace{-5pt}/\prod s_i[0]$.
The
corresponding field is defined, of course, as
$\prod_i Y(s_i,z)^{-1}$ in agreement with Proposition \ref{loc}.
Let us check vacuum axiom of vertex algebras for this field.
First of all, when we apply it to the vacuum, which has
$L[0]$ eigen-value zero in any graded vertex algebra, only non-positive
modes of $s_i$ survive. As a result, at $z=0$ we obtain
precisely $|0\vac\hspace{-5pt}/\prod s_i[0]$. Operator $T$ extends naturally,
because is commutes with $s[0]$. The second part of the vacuum axiom
also holds, because $$\{T,\prod_i s^{-1}(z)\}_-=
\{T,1/(\prod s_i)(z)\}_-=\prod s_i^{-1}(z)\{\prod s_i(z),T\}_-
\prod s_i^{-1}(z)
$$
$$
= - \prod s^{-2}(z)(\partial_z \prod s_i(z))=
\partial z \prod s_i^{-1}(z).
$$
In addition one can show that 
new fields are mutually local with all old fields 
$Y(b,z)$ and with each other. 

We will now construct a field for every element in $V_S$.
We will do it by induction on $l$ in filtration $F^lV$
of Proposition \ref{filtr}.
For an arbitrary element $a/\prod s_i[0]$ of grading $h$ we define
$$Y(a/\prod s_i[0],z)=\,
\no 
Y(|0\vac\hspace{-4pt}/\prod s_i[0],z)
Y(a,z)
\no
$$
$$
=
(\sum_{k\leq 0}(|0\vac\hspace{-5pt}/\prod s_i[0])[k]z^{-k})(\sum_i
a[i]z^{-i-h}) 
+
(\sum_i a[i]z^{-i-h}) (\sum_{k\geq 1}(|0\vac\hspace{-5pt}/\prod
s_i[0])[k]z^{-k}).
$$
One can see that this expression is well defined as an 
element in ${\rm End}V[[z,z^{-1}]]$. 
When we apply this field to the vacuum, second term does not
contribute, and we get $(\prod s_i^{-1})[0] a$. This is {\it not}\
the same as $a/\prod s_i[0]$ but the difference lies in the 
deeper part of the filtration.

So we have now constructed a field for each $a\in V_S$. To show
that this definition is compatible with changes from $a$ 
to $s[0]a/s[0]$, notice that all constructed fields are mutually
local and satisfy the second part of vacuum axiom. Then the argument 
of the uniqueness theorem of Section 4.4 of \cite{Kac} works and 
allows us to conclude that for any two ways of writing an element in $V_S$
the corresponding fields are the same. \hfill$\Box$

Above proposition allows us to define sheaves of 
vertex algebras over any algebraic variety $X$.
\begin{dfn}
{\rm
A (quasi-)loco sheaf ${\cal V}$ of vector spaces over $\CC$ is 
said to be a  
{\it a (quasi-)loco sheaf of vertex algebras}\ if for every affine
subset ${\rm
Spec}(R)\subset
X$ sections $\Gamma({\rm Spec}(R),{\cal V})$ form a
vertex $R$-algebra and restriction maps are precisely localization
maps of Proposition \ref{vertloc}.
}
\label{qlcvert}
\end{dfn}

Our goal now is to provide $H^*({\cal V})$ with a structure
of vertex algebra. Certainly, operators $T$ and $L[0]$, as well as 
the vacuum, behave well under localizations and are therefore globally
defined. For every affine set in $X$ and every integer $n$ we may consider
operation $(n)$ that maps $(a,b)$ to $a_{(n)}b$. It also commutes 
with localization, which gives us the map $(n):{\cal V}\otimes{\cal V}
\to {\cal V}$. This map induces a cup product on the cohomology
of ${\cal V}$ and we shall show that combining all these maps together
yields a vertex algebra structure on $H^*({\cal V})$. 
\begin{prop}
{\rm Cohomology of quasi-loco sheaf of vertex algebras ${\cal V}$ has a
natural structure of vertex algebra. Moreover, if sections of ${\cal V}$
over Zariski open sets are given the structure of conformal algebras 
that is compatible with localization maps, then
$H^*({\cal V})$ has a natural conformal structure.
}
\label{vertcoh}
\end{prop}

{\em Proof.} We will use the equivalent set of axioms of vertex algebras,
see \cite{Kac}, Section 4.8, This set consists of the partial vacuum
identity
$$Y(|0\vac,z)={\bf 1},~~a{(-1)}|0\vac=a$$ 
and Borcherds identity
$$\sum_{j=0}^\infty 
\left(\hspace{-4pt}
\begin{array}{c}  m \\ j
\end{array}
\hspace{-4pt}
\right)
(a_{(n+j)}b)_{(m+k-j)}c=
\sum_{j=0}^\infty (-1)^j
\left(\hspace{-4pt}
\begin{array}{c}  n \\ j
\end{array}
\hspace{-4pt}
\right)
a_{(m+n-j)}(b_{(k+j)}c)
$$
$$
-(-1)^{{\rm parity}(a){\rm parity}(b)}
\sum_{j=0}^\infty (-1)^{j+n}
\left(\hspace{-4pt}
\begin{array}{c}  n \\ j
\end{array}
\hspace{-4pt}
\right)
b_{(n+k-j)}(a_{(m+j)}c)
$$
for all $a,b,c\in V$ and all $k,m,n\in \ZZ$. When we have 
a graded algebra with $L[0]\geq 0$, for any given $L[0]$-eigen-values
of $a$, $b$ and $c$, the sums in Borcherds identity are finite.
Therefore, Borcherds identity is just a collection of identities
for maps $(n)$ between components ${\cal V}_r$ of ${\cal V}$.
Therefore, they induce identities on cohomology of ${\cal V}$
when we replace $(n)$ with corresponding cup products. 
A careful examination of signs shows that Borcherds identity holds
on cohomology of ${\cal V}$ if we define parity on $H^s({\cal V})$
as the sum of $s$ and parity on ${\cal V}$.
Partial vacuum identity on cohomology also follows from the vacuum
identity on ${\cal V}$ and the fact that $|0\vac$ is globally 
defined. If ${\cal V}$ has conformal structure, then $a=L[-2]|0\vac$
is globally defined, so it lies in $H^0$ and provides cohomology
with conformal structure of the same central charge.
\hfill$\Box$

We finish this section with discussion of {\it BRST
cohomology}\ . Let $V$ be a vertex algebra
and let $a$ be
an element of $V$ such that $a_{(0)}^2=0$.
Consider cohomology of $V$ with respect to $a_{(0)}$, which is called
BRST cohomology. Operator $a_{(0)}$ and field $Y(a,z)$
are called BRST operator and BRST field respectively.
The following proposition is standard.
\begin{prop}
{\rm 
BRST cohomology of $V$ with respect to $a_{(0)}$ has a
natural structure of vertex algebra.
}
\label{BRSTvert}
\end{prop}

{\em Proof.}
One has the following identity (\cite{Kac}, equation 4.6.9)
$$\{a_{0},Y(b,z)\}_\pm=Y(a_{(0)}b,z).$$
Therefore, if $b$ is annihilated by $a_{(0)}$ then
$Y(b,z)$ commutes with $a_{(0)}$ and conserves the kernel 
and the image of BRST operator. This provides us with the set
of mutually local fields for BRST cohomology, and it remains to
employ the uniqueness theorem of Section 4.4 of \cite{Kac}.
\hfill$\Box$

In particular, if $a$ is an odd element of $V$ such that all modes 
$a_{(n)}$ anticommute with each other, then $a_{(0)}$ could serve
as BRST operator. All major results of this paper involve BRST
cohomology by operators of this type.

\section{Chiral de Rham complex as a sheaf of vertex algebras over
a smooth variety}
In this section we review and summarize results of the 
extremely important paper of Malikov, Schechtman and 
Vaintrob \cite{MSV}. Notations of our paper follow closely
those of \cite{MSV}. We assume some familiarity with the vertex 
algebras of free bosons and fermions. The reader is referred to
\cite{Kac} or pretty much any conformal field theory textbook.

For every smooth variety $X$ authors of \cite{MSV} define a loco sheaf
of vertex algebras $\msv(X)$ which is called
{\it chiral de Rham complex} of $X$.
It is described in local coordinates $x^1,...,x^{\dim X}$
as follows.
There are given $2\dim X$ fermionic fields $\varphi^i(z)$, $\psi_i(z)$ 
and $2\dim X$ bosonic fields $a_i(z)$, $b^i(z)$ where index 
$i$ is allowed to run from $1$ to $\dim X$. The non-trivial
super-commutators between modes are given by 
$$\{a_i[k],b^j[l]\}_-=\delta_i^j\delta_{k+l}^0$$
$$\{\psi_i[k],\varphi^j[l]\}_+=\delta_i^j\delta_{k+l}^0,$$
and the fields are defined as 
$$a_i(z)=\sum_k a_i[k]z^{-k-1},~
b^i(z)=\sum_k b^i[k]z^{-k},$$
$$
\varphi^i(z)=\sum_k \varphi^i[k]z^{-k},~
\psi_i(z)=\sum_k \psi_i[k]z^{-k-1}.$$
There is defined a Fock space generated from the vacuum vector
$|0\vac$ by  non-positive modes of $b$ and $\varphi$ and by 
negative modes of $a$ and $\psi$. To obtain sections
of $\msv(X)$ over a neighborhood $U$ of $x=0$, one considers
tensor product $V$ of the above Fock space with the ring
of functions over $U$ with $b^i[0]$ plugged in place of $x^i$.
Tensor product is taken over the ring $\CC[b[0]]$. Grading on
this space is defined as the opposite of the sum of mode numbers.

Certainly one needs to specify how elements of $V$
change under a change of local coordinates.
For each new set of coordinates 
$$\tilde x^i=g^i(x),~x^j=f^j(\tilde x^j)$$
this is accomplished in \cite{MSV} by the formulas 
$$
\tilde b^i(z)=g^i(b(z))
$$
$$
\tilde\varphi^i(z)=g^i_j(b(z)) \varphi^j(z)
$$
$$\tilde a_i(z)=\,\no a_j(z) f^j_i(b(z))\no
+\no \psi_k(z) f^k_{i,l}(b(z))g^l_r(b(z))\varphi^r(z)\no
$$
$$
\tilde \psi_i(z)= \psi_j(z)f_i^j(b(z))
$$
where
$$g^i_j=\partial g^i/\partial x^j,
~f^i_j=(\partial f^i/\partial \tilde x^j)\circ g,
~f^i_{j,k}=(\partial^2 f^i/\partial \tilde x^j\partial \tilde x^k)
\circ g
$$
and normal ordering $\no\no$ is defined by pushing all positive 
modes of $a$, $b$, $\psi$ and $\varphi$ to the right, multiplying
by $(-1)$ whenever two fermionic modes are switched.

For any choice of local coordinates, one introduces fields
$$L(z)=\,\no \partial_z b^i(z)a_i(z)\no+\no\partial_z
\varphi^i(z)\psi_i(z)\no$$
$$J(z)=\,\no\varphi^i(z)\psi_i(z)\no$$
$$G(z)=\partial_z b^i(z)\psi_i(z)$$
$$Q(z)=\partial_z a_i(z)\varphi^i(z).$$
Field $L(z)$ is invariant under the change of coordinates, which 
provides $\msv(X)$ with structure of sheaf of conformal vertex
algebras. The $L[0]=0$ part is naturally isomorphic to the 
usual de Rham complex on $X$, with grading given by 
eigen-values of $J[0]$ and differential given by $Q[0]$ (both modes 
are globally defined). If $X$ is Calabi-Yau, all four of the above
fields are well-defined which provides $\msv(X)$ with structure
of sheaf of {\it topological vertex algebras.} 
This means that spaces of sections over affine subsets are equipped with
structures of topological vertex algebras, in a manner consistent
with localization. This structure is analogous to the conformal
structure but requires a choice of four fields $Q$, $G$, $J$ and $L$
that satisfy certain OPEs, see \cite{MSV}.

It was suggested in \cite{MSV} that cohomology of
$\msv(X)$ has a structure of vertex algebra which describes
the holomorphic part of A model of $X$, see \cite{witten}.
Since we now know how to provide cohomology of the loco sheaf
of vertex algebras $\msv(X)$ with such structure, we can
state the following definition.
\begin{dfn}
{\rm
Let $X$ be a smooth algebraic variety over $\CC$. We define
{\it A model topological vertex algebra}\ of $X$ to be $H^*(\msv(X))$
with the structure of vertex algebra on it defined in
Proposition \ref{vertcoh}. This algebra also possesses
the conformal structure since $L(z)$ is globally defined,
as well as the structure of the topological algebra,
with operators given by formulas of \cite{MSV} in the case when
$X$ is Calabi-Yau.
}
\label{amodel}
\end{dfn}

If $X$ is a Calabi-Yau variety, one can also talk about B model
topological vertex algebra of $X$. As a vertex algebra, it is identical to
the A model, but
the additional structure of the topological algebra differs. 
\begin{dfn}
{\rm
Let $X$ be a smooth Calabi-Yau manifold. We define 
{\it B model topological vertex algebra}\ of $X$ as follows. As a vector
space, it coincides with the A model topological vertex algebra of $X$.
Operator $T$ is also the same, and so are the fields of 
the algebra. 
Topological structure of the B model vertex algebra 
is related to the topological structure of A model algebra
by {\it mirror involution} 
$$Q_B=G_A,~G_B=Q_A,~J_B=-J_A,~L_B=L_A-\partial J_A.$$
}
\label{bmodel}
\end{dfn}
In what follows we will often abuse the notations and call these 
algebras simply A and B model vertex algebras, but it should always
be understood that they are considered together with their extra
structure. Notice that B model is ill-defined for varieties $X$ that are
not Calabi-Yau, even as a conformal vertex algebra, because $L_B$ is
ill-defined for them.

The first goal of our paper is to calculate A  and B model
vertex algebras for Calabi-Yau hypersurfaces in smooth
toric nef-Fano varieties. The second goal (which is
only partially achieved) is to generalize our
results to toric varieties with Gorenstein singularities. 

\section{Vertex algebras of line bundles and zeros of their
sections}

In this section we study \MSV\ of a line bundle $L$
over a smooth variety $P$. Given a section of the dual line
bundle we are able to calculate the \MSV\ of its zero set in
terms of the push-forward of the \MSV\ of the line bundle.

We denote the projection to the base
by $\pi: L\to P$. The
line bundle structure is locally
described by the fact that there is one special coordinate
$x^1$ such that allowed local changes of coordinates
are compositions of changes in $x^2,...,x^{\dim L}$ and
changes 
$$\tilde x^1\,=\,x^1\, h(x^2,...,x^{\dim L}),~~
\tilde x^i=x^i,~i\geq 2.$$

%Conservation of B_0\PSI_0 (z)
\begin{prop}
{\rm Field $b^1(z)\psi_1(z)$ depends only on the line bundle
structure of $L$.
}
\label{beepsi}
\end{prop}
{\em Proof.} This field is clearly unaffected by any changes
of coordinates on the base that leave $x^1$ intact. In addition,
for the change of coordinates as above, we have
$$\tilde b^1=b^1 h(b^2,...,b^{\dim L}),~
\tilde \psi_1=\psi_1/h(b^2,...,b^{\dim L}).$$
As a result, field $b^1(z)\psi_1(z)$ is independent from 
the choice of local coordinates that are compatible with the given 
line bundle structure. \hfill $\Box$

%CY condition
\smallskip
The following lemma is clear.
\begin{lem}
{\rm Line bundle $L$ has trivial canonical class if and only if
$L$ is the canonical line bundle on $P$.
}
\label{CYlb}
\end{lem}

%remark on push-forwards
\begin{rem}
{\rm Since $\msv(L)$ is a loco sheaf, its 
cohomology spaces $H^* (\msv(L))$ are isomorphic to cohomology
$H^*(\pi_*\msv(L))$ of its push-forward to $P$, which is the sheaf 
we are mostly interested in. 
}
\label{pushforward}
\end{rem}

%remark on odd bundles
\begin{rem}
{\rm One may also consider the bundle $\Pi L$ obtained by declaring
the coordinate on $L$ odd. It turns out that the corresponding sheaf of
vertex algebras is roughly the same as the corresponding
sheaf for the even bundle $L^{-1}$. More precisely, the push-forwards of both
bundles to $P$ coincide. Locally the isomorphism is obtained by mapping
$(b^1,\varphi^1,a_1,\psi_1)$ for $\Pi L$ to
$(\psi_1,a_1,\varphi^1,b^1)$ for $L^{-1}$. 
It has been observed in
\cite{Schwarz} that mirror symmetry for Calabi-Yau complete intersections
may be formulated in terms of odd bundles on ambient projective
varieties.
}
\end{rem}

Let us additionally assume that we have in our disposal a section 
$\mu$ of the dual line bundle $L^{-1}$. This amounts to having a global
function on $L$ which is linear on fibers. We will also assume 
that zeros of $\mu$ form a reduced non-singular divisor $X$ on
$P$. The goal of the rest of the section is to describe 
$\msv(X)$ in terms of $\pi_*\msv(L)$.

%conservation of d(section)(z)
The following lemma is easily checked by a calculation in local
coordinates. In fact, it holds for any global function on any smooth
variety.
\begin{lem}
{\rm 
Fields $\mu(z)$ and $D\mu(z)$ that are locally defined as
$$\mu(z)=\mu(b^1,...,b^{\dim L})(z),~~
D\mu(z)=\sum_i \varphi^i(z) {{\partial \mu}\over{\partial b^i}}(z)
$$
are independent of the choice of coordinates and are therefore globally 
defined. In particular, the operator ${\cal BRST}_\mu=\oint D\mu(z)dz$ is
globally
defined.
}
\label{dmu}
\end{lem}

%BRST-type description of the VOA of the zeros of the section
It is now time to state the main result of this section.
\begin{theo}
{\rm 
Let $X$ be a smooth hypersurface in a smooth variety $P$ 
defined as above by a section $\mu$ of line bundle $L^{-1}$.
Then sheaf of vertex algebras $\msv(X)$ is isomorphic to BRST
cohomology of sheaf $\pi_*\msv(L)$ with respect to
the operator ${\cal BRST}_\mu$. Here cohomology is
understood in the sense of sheaves, that is as a 
sheafification of BRST cohomology presheaf.
}
\label{BRSThyper} 
\end{theo}

{\em Proof.} Clearly, ${\cal BRST}_\mu$ is a differential,
because its anticommutator with itself is zero. We denote
the BRST cohomology of $\pi_*\msv(L)$ by $\widehat{\msv(X)}$ and it
is our goal to construct an isomorphism between $\widehat{\msv(X)}$
and $\msv(X)$. 

It is enough to construct this isomorphism locally
for any point $p\in P$ provided that our construction withstands
a change of coordinates. We use Hausdorff
topology on $P$ rather than Zariski topology. Point $p$ 
may or may not lie in $X$ so our discussion splits into two
cases.

{\em Case 1.} $p\notin X$

In this case for any sufficiently small neighborhood $U\subset P$ of $p$
we can choose a coordinate system $(x^1, x^2,...,x^{\dim L})$
on $\pi^{-1}U$ such that $x_1$ is the special line bundle variable,
and $\mu=x^1$. As a result, 
$${\cal BRST}_\mu=\oint 
\varphi^1(z) dz=\varphi^1[1].$$
A simple calculation on the flat space then shows that 
cohomology by ${\cal BRST}_\mu$ are zero on 
sections of $\pi_*\msv(L)$ for any sufficiently small $U$.
As a result, $\widehat{\msv(X)}$ is supported on $X$, which
is, of course, true of $\msv(X)$. 

{\em Case 2.} $p\in X$. 

For any sufficiently small neighborhood $U\subset P$ of $p$ we can 
choose a system of coordinates $(x^1,x^2,...,x^{\dim L})$ 
on $\pi^{-1}U$ that agrees with line bundle structure, such that
$\mu=x^1x^2$ and $x^3,...,x^{\dim L}$ form a system of local
coordinates on $X\cap U$. We then have
$${\cal BRST}_\mu = \oint (b^1 \varphi^2 + b^2 \varphi^1)(z)dz
=\sum_{k\in\ZZ} (b^1[k]\varphi^2[-k+1]+b^2[k]\varphi^1[-k+1]).$$
Fock space $\Gamma(U,\pi_*\msv(L))$ is a tensor product
of spaces ${\Fock}_{1,2}$ and ${\Fock}_{\geq3}$,
which are the spaces generated by modes of $a_i,b^i,\varphi^i,\psi_i$ 
for $i\in\{1,2\}$ and $i\in\{3,...,\dim L\}$ respectively.
Since  ${\cal BRST}_\mu$ acts on the first component
of this tensor product, its cohomology is isomorphic to the tensor product 
of ${\Fock}_{\geq3}$ and cohomology of ${\Fock}_{1,2}$
with respect to ${\cal BRST}_\mu$.

We claim that cohomology of  ${\Fock}_{1,2}$
with respect to ${\cal BRST}_\mu$ is one-dimensional and is generated
by the image of vacuum vector $|0\vac$. We do not multiply by 
$\Gamma(U,{\cal O}_P)$, which does not alter the argument.
Notice first that $\Fock_{1,2}$ is a restricted tensor product (that
is almost all factors
are $1$) of the following infinite set of vector spaces.

$$\bullet\oplus_{l\geq 0} \CC(a^1[-k])^l
\oplus_{l\geq 0} \CC(a^1[-k])^l\varphi^2[-k+1],
~~{\rm for~all}~ k>0$$

$$\bullet\oplus_{l\geq 0} \CC(a^2[-k])^l
\oplus_{l\geq 0} \CC(a^2[-k])^l\varphi^1[-k+1],
~~{\rm for~all}~ k> 0$$

$$\bullet\oplus_{l\geq 0} \CC(b^1[-k])^l
\oplus_{l\geq 0} \CC(b^1[-k])^l\psi^2[-k-1],
~~{\rm for~all}~ k\geq 0$$

$$\bullet\oplus_{l\geq 0} \CC(b^2[-k])^l
\oplus_{l\geq 0} \CC(b^2[-k])^l\psi^1[-k-1],
~~{\rm for~all}~ k> 0$$

$$\bullet\oplus_{l\geq 0} R(b^2[-k])^l
\oplus_{l\geq 0} R(b^2[-k])^l\psi^1[-k-1],
~~{\rm for~all}~ k> 0$$

In the last formula $R$ means the ring of function on a disc. We assume
here that the neighborhood $U$ is a product of $|x^2|\leq c$ and 
some $U_{x^3,...,x^{\dim L}}$.

Vacuum vector, of course, corresponds to the product of all $1$.
The Fock space is graded by the eigen-values of $L_{1,2}[0]$
that is by the opposite of the total sum of indices. 
Operator ${\cal BRST}_\mu$ shifts this grading by $-1$.
If we consider elements with bounded grading, it is enough
to consider only product of a finite number of above spaces.
For each such product, ${\cal BRST}_\mu$ is a sum of anticommuting
operators on each component. One can then show that cohomology is a
tensor product of cohomologies for each component by induction 
on the number of components. On each step of the induction we
are using a spectral sequence for a stupid filtration of
the tensor product complex, with grading given by eigen-values of 
$L_{1,2}[0]$.

As a result, to show that cohomology space is one-dimensional,
it is enough to show that for each of the spaces above cohomology
is one-dimensional and is given by the image of $1$.
It is sufficient to consider the first, third, and fifth types only.
For a space of first type, kernel of ${\cal BRST}_\mu$ is
$$\CC 1\oplus~\oplus_{l\geq 0} \CC(a^1[-k])^l\varphi^2[-k+1]$$
and its image is 
$$\oplus_{l\geq 0} \CC(a^1[-k])^l\varphi^2[-k+1]$$
so the image of $1$ generates cohomology.
For a space of third type,  kernel is 
$$\oplus_{l\geq 0} \CC(b^1[-k])^l$$
and image is 
$$\oplus_{l\geq 1} \CC(b^1[-k])^l$$
which gives the same result.
For the space of the fifth type, we use $R/xR=\CC$.

So we managed to show that for a given choice of coordinates
on $\pi^{-1}U$, there is an isomorphism between sections
of $\widehat{\msv(X)}$ and $\msv(X)$. The proof is not over yet, 
because we need to show that these locally defined isomorphisms
could be glued together. This amounts to the demonstration that
the isomorphism just constructed commutes with any changes of
coordinates on $\pi^{-1}U$ that preserve our setup. Every such coordinate
change could be written in the form
$$
\tilde x^1=x^1\cdot h(x^2,...x^{\dim L}),~
\tilde x^2=x^1/h(x^2,...,x^{dim L})
$$
$$
\tilde x^i=f^i(x^3,...,x^{\dim L})+x^2g^i(x^2,...,x^{\dim L}),~
i\geq 3.
$$
It is clear that when $h=1$ and $g^i=0$ the corresponding splitting of
the Fock space is unaffected and the resulting isomorphism precisely
matches the change of variables on $X$. As a result, we only need to
show that isomorphism commutes with coordinate changes such
that $f^i(x)=x^i$. 

One can show that in this case fields $\tilde a_i(z),\tilde b^i(z),
\tilde \varphi^i(z),\tilde \psi_i(z)$ for $i\geq 3$ act on the cohomology
in the same way as the operators $a_i(z),b^i(z),\varphi^i(z),\psi_i(z)$,
because the difference lies in the image of ${\cal BRST}_\mu$. This
finishes the proof. \hfill$\Box$

It is clear that our isomorphism commutes with structures of
sheaves of vertex algebras.
We also have the following corollary which will be very useful later.
\begin{prop}
{\rm 
For any affine subset $U\subset P$ the ${\cal BRST}_\mu$ cohomology
space of $\Gamma(U,\pi_*\msv(L))$ is isomorphic to $\Gamma(U,\msv(X))$.
}
\label{affine}
\end{prop}

{\em Proof.} Sheaf $\pi_*\msv(L)$ is a quasi-loco sheaf, and
${\cal BRST}_\mu$ is a map of a quasi-loco sheaf into itself.
For any affine subset it is induced from the map of corresponding
loop-modules, and then everything follows from the Remark \ref{msh}.
\hfill $\Box$

%calculation of the N=2 super-algebra for the CY case
We are especially interested in the case where $L$ has a non-degenerate
top form. In this case, by Lemma \ref{CYlb}, $L$ is a
canonical line bundle, and a section $\mu$
of $L^{-1}$ produces a Calabi-Yau divisor $X$ on $P$.
Our goal here is to calculate global fields $G(z)$ and $Q(z)$
on $\msv(X)$ in terms of some global fields on $\pi_*\msv(L)$.

\begin{prop}
{\rm When $L$ is canonical bundle on $X$,  field $G_X(z)$ is the
image of field $G_L(z)-(b_1(z)\psi_1(z))'$. Field $Q_X(z)$ 
is the image of field $Q_L(z)$.
}
\label{GZhyper}
\end{prop}

{\em Proof.} Because of \ref{beepsi}, all fields in question are
defined globally, so a local calculation is sufficient. We assume
notations of the proof of Theorem \ref{BRSThyper}. Then we have 
$$Q_X(z)-Q_L(z)=-a_1(z)\varphi^1(z)-a_2(z)\varphi^2(z),
$$$$
G_X(z)-G_L(z)+(b^1(z)\psi_1(z))'=b^{2}(z)'\psi_2(z)-b^1(z)\psi_1(z)',$$
and we need to show that right-hand sides of these equations
are commutators of ${\cal BRST}_\mu$ and some fields.
This goal is accomplished by fields
$-a_1(z)a_2(z)$ and $-\psi_1'(z)\psi_2(z)$ respectively.
\hfill$\Box$

\section{BRST description of vertex algebra in logarithmic coordinates}
This section is in a sense a mirror of the previous one. It contains
a local calculation of \MSV\ of a smooth toric variety as BRST cohomology
of some $\cal MSV$-like space defined in terms of local coordinates.

We introduce some notations which will stay with us for
the rest of the paper. Let $M$ be a free abelian group of rank
$\dim M$ and $N={\rm Hom}(M,\ZZ)$ be its dual. 
Vector space $(M\oplus N)\otimes \CC$ has dimension $2\dim M$ and it
is equipped with a standard bilinear form denoted
by"$\cdot$". This allows us to construct $2\dim M$ bosonic and $2\dim M$
fermionic fields. Really, one can always construct $k$ bosonic and 
$k$ fermionic fields starting from a vector space or dimension $k$
with a non-degenerate bilinear form on it, see for example \cite{Kac}, 
so our purpose here is to fix notations. 
For every $m\in M$ and $n\in N$ we have 
$$m\cdot B(z)=\sum_{k\in \ZZ}m\cdot B[k]z^{-k-1},
~n\cdot A(z)=\sum_{k\in \ZZ}n\cdot A[k]z^{-k-1},$$
$$m\cdot \Phi(z)=\sum_{k\in \ZZ}m\cdot \Phi[k]z^{-k},
~n\cdot \Psi(z)=\sum_{k\in \ZZ}n\cdot \Psi[k]z^{-k-1}.$$
Notice that the moding of $B$ also has $z^{-k-1}$ in it, in contrast to
the moding of $b^i$ in the previous section.
The non-zero super-commutators are 
$$\{m\cdot B[k],n\cdot A[l]\}_- = (m\cdot n) k\delta_{k+l}^0 {\bf id},$$
$$\{m\cdot \Phi[k],n\cdot \Psi[l]\}_+ = 
(m\cdot n) \delta_{k+l}^0 {\bf id}.$$

Our battlefield will be the following space whose construction
is standard as well.
$$\Fock_{M\oplus N}=_{{\rm def}} \oplus_{m\in M,n\in N}
~\otimes_{k\geq 1} \CC[B[-k]] 
~\otimes_{k\geq 1} \CC[A[-k]] 
~\otimes_{l\geq 0} (\CC+\CC\Phi[-l]) 
$$
$$
~\otimes_{l\geq 1} (\CC+\CC\Psi[-l]) |m,n\vac$$
Here $\otimes$ means restricted tensor product over $\CC$ where
only finitely many factors are not equal to $1$. Vectors 
$|m,n\vac$ are annihilated by positive modes of $A$, $B$, and $\Phi$,
and by non-negative modes of $\Psi$. Also, 
$$A[0]|m,n\vac=m|m,n\vac,~B[0]|m,n\vac=n|m,n\vac.$$

This Fock space possesses a structure of vertex algebra,
see for example \cite{Kac}. Among the fields
of this 
algebra the important role is played by so-called {\it vertex operators}\
$$\no\ee^{\int (m\cdot B(z)+n\cdot A(z))}\no$$
which are defined as follows 
$$\no\ee^{\int (m\cdot B(z)+n\cdot A(z))}\no
\prod A[...]\prod B[...]\prod
\Phi[...]\prod\Psi[...]|m_1,n_1\vac
$$
$$
=C(m,n,m_1,n_1)
z^{m\cdot n_1+n\cdot m_1}
\prod_{n<0} \ee^{-(m\cdot B[n]+n\cdot A[n])
{z^n\over n}} 
\prod_{n>0} \ee^{-(m\cdot B[n]+n\cdot A[n])
{z^n\over n}}
\prod A[...]
$$
$$\prod B[...]\prod
\Phi[...]\prod\Psi[...]
|m+m_1,n+n_1\vac.$$
Cocycle $C(m,n,m_1,n_1)$ here equals 
$(-1)^{m\cdot n_1}$. It is used to make vertex operators
purely bosonic. Our notation 
suppresses this cocycle, which is a bit unusual 
but should not lead to any confusion.
Vertex operators obey the following OPEs
$$\no\ee^{\int (m\cdot B(z)+n\cdot A(z))}\no
\no\ee^{\int (m_1\cdot B(w)+n_1\cdot A(w))}\no
={\no\ee^{\int (m\cdot B(z)+n\cdot A(z))}
\ee^{\int (m_1\cdot B(w)+n_1\cdot A(w))}
\no\over (z-w)^{m\cdot
n_1+n\cdot m_1}}$$
where putting both fields under the same $\no\no$ sign
means that we move all negative modes to the left and
all positive modes to the right as in the definition of 
vertex operators above. Of course, this OPE could be expanded 
by Taylor formula, and the resulting fields are normal
ordered products of vertex operators, free bosons, and
their derivatives. In general, all fields of the vertex algebra
$\Fock_{M\oplus N}$ are normal ordered products of various
$B$, $A$, $\Psi$, $\Phi$ and their derivatives times one (perhaps
trivial) vertex operator.

\begin{rem}
{\rm
Vertex algebra $\Fock_{M\oplus N}$ possesses a conformal
structure, given by 
$$L_{M\oplus N}(z)=\,\no B(z)\cdot A(z)\no
+\no(\partial_z\Phi(z)\cdot\Psi(z)\no.$$
The corresponding grading operator $L_{M\oplus N}[0]$ assigns grading $m\cdot n$ to
vector $|m,n\vac$. Elements $\prod_im_i\cdot\Phi[0]|m,n\vac$ have the
same eigenvalue, but for every other element with the same 
$A[0]$ and $B[0]$ eigen-values, the grading is strictly larger
than $m\cdot n$. Under this conformal structure the moding of the fields
$A$, $B$, $\Psi$ and $\Phi$ is as above.
}
\label{lflat}
\end{rem}

%dimension one case
The goal of this section is to construct a flat space vertex algebra
in terms of $A,B,\Phi,\Psi$.  We first look at the case of dimension
one. In this case $M$ and $N$ are one-dimensional, and
$A,B,\Phi,\Psi$ are no longer vector-valued. Consider the 
following fields 
$$b(z)=\ee^{\int B(z)},~\varphi(z)=\Phi(z)\ee^{\int
B(z)}, ~\psi(z)=\Psi(z)\ee^{-\int B(z)},$$
$$a(z)=\,\no A(z)\ee^{-\int B(z)}\no +\no \Phi(z)\Psi(z)\ee^{-\int
B(z)}\no.$$ 

\begin{prop}
{\rm 
Operator product expansions of $a,b,\varphi,\psi$ are 
$$a(z)b(w)={1\over{z-w}}+\reg,~~\varphi(z)\psi(w)=
{1\over{z-w}}+\reg,$$
and all other OPEs are non-singular.
}
\label{expvoa}
\end{prop}

{\em Proof.} It is a standard calculation of OPEs that include
vertex operators, which we omit. This is not too surprising since
the fields in question are given by the formulas of \cite{MSV} applied 
to the exponential change of variables $\tilde x={\rm exp}(x)$.
\hfill$\Box$

\begin{prop}
{\rm Modes of the fields $a,b,\varphi,\psi$ generate a vertex algebra
which is isomorphic to global sections of \MSV\ of a one-dimensional
affine space.
}
\label{smallfock}
\end{prop}

{\em Proof.} First, all OPEs are correct due to \ref{expvoa}.
Since $a$ and $b$ are bosonic and $\varphi$ and $\psi$ are fermionic,
this implies that super-commutators of their modes are correct.
Notice that the conformal weight of $b$ is zero, so it is moded correctly.
One can also show that positive modes of $b,\varphi$ and
non-negative modes of $a,\psi$ annihilate $|0,0\vac$. The rest follows
from the fact that Fock representation of the algebra of modes is
irreducible. \hfill$\Box$

The following calculation is extremely useful.
\begin{prop}
{\rm We define $L(z)$, $J(z)$, $Q(z)$, and $G(z)$ for
$a,b,\varphi,\psi$ as usual, see \cite{MSV}. Then in terms
of $A,B,\Phi,\Psi$ we have
$$
Q(z)=A(z)\Phi(z)-\partial_z\Phi(z),~G(z)=B(z)\Psi(z),
$$
$$
J(z)=\,\no\Phi(z)\Psi(z)\no + B(z)
$$
$$L(z)=\,\no B(z)A(z)\no +\no \partial_z\Phi(z)
\Psi(z)\no.$$\\[-3em]
}
\label{toricVOA}
\end{prop}

{\em Proof.} It is a standard calculation, which is again omitted. 
\hfill$\Box$

We now define $\Fock_{M\oplus N_{\geq 0}}$ as a subalgebra of 
$\Fock_{M\oplus N}$ characterized by the condition that 
eigen-values of $B[0]$ are non-negative. This amounts to only
allowing $|m,n\vac$ with $n\geq 0$.
We will now show that vertex algebra generated by $a,b,
\varphi,\psi$ could be obtained as a certain 
BRST cohomology of vertex
algebra $\Fock_{M\oplus N_{\geq 0}}$.

\begin{theo}
{\rm
Vertex algebra of $a,b,\varphi,\psi$ is isomorphic to BRST cohomology
of $\Fock_{M\oplus N_{\geq 0}}$ with respect to the operator
$${\cal BRST}_g=\oint {\cal BRST}_g(z)dz=\oint g \Psi(z)\ee^{\oint
A(z)}dz$$
where $g$ is an arbitrary non-zero complex number.
}
\label{dimone}
\end{theo}

{\em Proof.}
First of all, notice that all modes of $a,b,\varphi,\psi$
commute with ${\cal BRST}_g$. Really, all these fields except $a(z)$ give
non-singular OPEs with ${\cal BRST}_g(w)$, and
$$a(z){\cal BRST}_g(w)= {{g\no A(z)\Psi(w)\ee^{\int
-B(z)+A(w)}\no}\over{z-w}}+\reg
$$
$$+g(-{{\Psi(z)}\over{z-w}}+O(z-w)){{\ee^{-\int B(z)+\int A(w)}}\over
{z-w}}
$$
$$
={{-g\Psi(z)\no\ee^{\int(A(z)-B(z))}\no}\over{(z-w)^2}}+\reg
$$
which implies $\{a(z),{\cal BRST}_g\}_-=0$.

Space $\Fock_{M\oplus N_{\geq 0}}$ is graded by eigen-values
of $B[0]$ and ${\cal BRST}_g$ shifts them by one. We
first show that ${\cal BRST}_g$ has no cohomology for eigen-values
of $B[0]$ that are positive. Really, we can look at the operator
${\cal R}(z)=\Phi(z)\ee^{-\int A(z)}$. A similar calculation
shows that 
$$\{{\cal R}(z),{\cal BRST}_g\}_+={\bf id}$$
and therefore the anticommutator of the zeroth mode of $R(z)$ 
and ${\cal BRST}_g$ is identity. Thus we found a homotopy
operator, which insures that there is no cohomology at positive
eigen-values of $B[0]$.

Fortunately, the above operator shoots out of $\Fock_{M\oplus N_{\geq 0}}$
from zero eigen-space of $B[0]$. So we found that the cohomology
is isomorphic to the kernel of ${\cal BRST}_g$ on the zero eigen-space
of $B[0]$. To show that all elements of this space can be obtained
by applying modes of $a,b,\varphi,\psi$ to $|0,0\vac$, we employ
the result of Proposition \ref{toricVOA}. More precisely,
$L[0]$ has non-negative eigen-values. Moreover, its zero eigen-space
is 
$$\oplus_{m\in\ZZ}(\CC\oplus\CC\Phi[0])|m,0\vac.$$
Since $L[0]$ commutes with ${\cal BRST}_g$, 
it induces grading on the kernel.

We prove by induction on eigen-values of $L[0]$
that all elements of the kernel of ${\cal BRST}_g$ with zero eigen-value
of
$B[0]$ are obtained by applying modes of $a,b,\varphi,\psi$ to
$|0,0\vac$. For $L[0]=0$ notice that cohomology is graded
by eigen-values of $A[0]$. An explicit calculation then shows that
for $k<0$ elements ${\cal BRST}_g |k,0\vac$ and ${\cal BRST}_g
\Phi[0] |k,0\vac$ are linearly independent. In addition,
${\cal BRST}_g \Phi[0]|0,0\vac$ is non-zero (it is proportional to
$|0,1\vac$),
and the rest is generated by modes of $b$ and $\varphi$.

If $L[0]v=lv$ with $l>0$, notice that 
$$L[0]=\sum_{k<0} ka[k]b[-k] + \sum_{k>0}kb[-k]a[k]
-\sum_{k<0} k\psi[k]\varphi[-k] + \sum_{k>0}k\varphi[-k]\psi[k].
$$
When applied to $v$, only finitely many terms survive.
So we have
$$v={1\over l}\sum_i p_iq_i v.$$
Since $v$ is in the kernel of ${\cal BRST}_g$, and $p_i$
commutes with ${\cal BRST}_g$, $p_i v$ is in the kernel for each $i$. 
Also $p_i v $ has strictly lower eigen-value of $L[0]$, so it is generated
by  modes of $a,b,\varphi,\psi$ due to the induction assumption.
Therefore, $v$ is also generated by modes of $a,b,\varphi,\psi$,
which finishes the proof.
\hfill$\Box$

%case of many dimensions
We can extend this theorem to lattices of any dimension as follows.
Consider a primitive cone $K^*$ in lattice $N$. Primitive here
means that it is generated by a basis $n_1,...,n_{\dim N}$ 
of $N$. The dual basis is denoted by $m_1,...,m_{\dim M}$.
We denote by $\Fock_{M\oplus K^*}$ the subalgebra of $\Fock_{M\oplus N}$
where eigen-values of $B[0]$ are allowed to lie in $K^*$. We consider
vertex algebra of flat space that is generated by fields
$$b^i(z)=\ee^{\int m_i\cdot B(z)},
~\varphi^i(z)=(m_i\cdot \Phi(z))\ee^{\int
m_i\cdot B(z)}, ~\psi_i(z)=(n_i\cdot \Psi(z))\ee^{-\int m_i\cdot B(z)},$$
$$a_i(z)=\,\no (n_i\cdot A(z))\ee^{-\int m_i\cdot B(z)}\no +\no
(m_i\cdot\Phi(z))(n_i\cdot\Psi(z))\ee^{-\int
m_i\cdot B(z)}\no$$ 
for all $i=1,...,\dim M$.

\begin{theo}
{\rm
Vertex algebra of $a_i,b^i,\varphi^i,\psi_i$ is isomorphic to BRST
cohomology
of $\Fock_{M\oplus K^*}$ with respect to operator
$${\cal BRST}_g=
\oint {\cal BRST}(z)dz=
\oint \sum_i g_i (n_i\cdot\Psi(z))\ee^{\oint n_i\cdot A(z)}dz$$
where $g_1,...g_{\dim M}$ are arbitrary non-zero complex numbers.
Moreover, operators $L$, $J$, $G$ and $Q$ are given by
$$
Q(z)=A(z)\cdot\Phi(z)-{\rm deg}\cdot \partial_z\Phi(z),~G(z)=
B(z)\cdot\Psi(z),
$$
$$
J(z)=\,\no\Phi(z)\cdot\Psi(z)\no +\, {\rm deg}\cdot B(z)
$$
$$L(z)=\,\no B(z)\cdot A(z)\no +\no \partial_z\Phi(z)\cdot
\Psi(z)\no$$
where "deg" is an element in $M$ that equals $1$ on all generators
of $K^*$.
}
\label{dimany}
\end{theo}

{\em Proof.} The result follows immediately from Theorem \ref{dimone}. 
Really, $\Fock_{M\oplus K^*}$ 
is a tensor product of $\dim M$ spaces discussed in there.
We grade each space by eigen-values of $m_i\cdot B[0]$,
and ${\cal BRST}_{g_i}$ becomes a degree one differential.
Operator ${\cal BRST}_g$ is a total differential on the corresponding
total complex, which finishes the proof. Another option is to
go through the proof of Theorem \ref{dimone} with minor changes
due to higher dimension. \hfill $\Box$

\section{Smooth toric varieties and hypersurfaces}
The result of the previous section could be interpreted as a calculation
of \MSV\  for a smooth affine toric variety given
by cone $K^*\subset N$. The first objective of this section is
to learn how to glue these objects together to get \MSV\ of 
a smooth toric variety (or rather a canonical line bundle over it). Then 
we employ the result of Theorem \ref{BRSThyper} to calculate vertex
algebras of Calabi-Yau hypersurfaces in toric varieties.

Let us recall the set of data that defines a smooth toric variety.
For general theory of toric varieties see \cite{danilov,fulton,oda}.
Paper of Batyrev \cite{bat.dual} may also be helpful.
A toric variety $\PP_\Sigma$ is given by {\it fan}\ $\Sigma$
which is a collection of rational polyhedral cones in $N$ with
vertex at $0$  such that 

(a) for any cone $C^*$ of $\Sigma$, $C^*\cap (-C^*)=\{0\}$;

(b) if two cones intersect, their intersection is a face in
both of them;

(c) if a cone lies in $\Sigma$, then all its faces lie in $\Sigma$,
this also includes the vertex (zero-dimensional face) of the
cone.

Toric variety is smooth if and only if all cones in $\Sigma$
are basic, that is generated by a part of a basis of $N$.

For every cone $C^*\in \Sigma$, one considers the dual cone $C\in M$
defined by
$$C=\{m\in M,~{\rm s.t.}~ m\cdot C^*\geq 0\}$$
and the corresponding affine variety $\AA_C={\rm Spec}(\CC[C])$.
We employ a multiplicative notation and denote elements
of $\CC[M]$ by $x^m$ for all $m\in M$. If $C_1^*$ is a face of $C_2^*$,
then $\CC[C_1]$ is a localization of $\CC[C_2]$ by all $x^m$ 
for which $m\in C_2, m\cdot C_1^*=0$. This allows us to construct
inclusion
maps between affine varieties $\AA_{C_i}$ and then glue them all
together to form a toric variety $\PP=\PP_\Sigma$.

We already know how to describe sections of \MSV\ on a flat space 
that corresponds to cone $K^*$ of maximum dimension. Our next step
is to construct vertex algebra that corresponds to a face of such cone.
Namely, let $C_1^*$ be generated by $n_1,...,n_r$ where 
$n_1,...,n_{\dim N}$ form a basis of $N$ and generate $C^*$. 
Then we can consider BRST operator 
$${\cal BRST}_g=\oint \sum_{i=1}^r g_i (n_i\cdot \Psi(z))
\ee^{\int n_i\cdot A_i(z)} dz$$
that acts on $\Fock_{M\oplus C_1^*}$. 
Corresponding BRST cohomology will be denoted by $\VOA_{C_1,g}$
and sections of \MSV\ on $\AA_C$ will be denoted by
$\VOA_{C,g}$. We have a natural surjective map
$$\rho : \Fock_{M\oplus C^*}\to \Fock_{M\oplus C_1^*}$$
which commutes with ${\cal BRST}_g$. Here we, of course,
abuse the notation a little bit by using two different definitions
of ${\cal BRST}_g$ for $C^*$ and $C_1^*$. However, we assume
that $g_i$ there are the same for $i=1,...,r$.

For every $m\in C$ and every $l\in \ZZ$ element $\ee^{\int m\cdot
B}[l]$ acts on both $\VOA_{C,g}$ and $\VOA_{C_1,g}$.
Really, its action on $\Fock_{M\oplus C^*}$ commutes with ${\cal BRST}_g$,
because $m\cdot C^*\geq0$. Consider multiplicative system $S$ 
generated by  elements $\ee^{\int m\cdot B(z)}[0]$ with $m\in
C,~m\cdot C_1^*=0$.

\begin{prop}
{\rm 
Map $\rho$ induces map 
$$\rho_{BRST}: \VOA_{C,g}\to\VOA_{C_1,g}$$
which is precisely the localization map of $\CC[C]$-loop-module
$\VOA_{C,g}$ with respect to multiplicative system $S$.
}
\label{localization}
\end{prop}

{\em Proof.} First we show that this map is the localization map
of corresponding vector spaces with the action of the multiplicative
system. For this it is enough to show that $\rho$ is the localization map.
This amounts to showing that any element $v$ of $\Fock_{M\oplus C}$ with
eigen-values of $B[0]$ equal to $n$ where $n\notin C_1^*$ is annihilated
by
some element in $S$. Since $n\notin C_1$, there exists an element $m\in
C$ such that $m\cdot C_1=0$ and $m\cdot n > 0$. It is easy to see that
a power of $s[0]= \ee^{\int m\cdot B}[0]$ annihilates $v$. Really, it
does not change the $L[0]$ eigen-value of $v$ but on the other hand it
increases its $A[0]\cdot B[0]$ eigen-value by an arbitrary positive
multiple 
of $m\cdot n$. So for big $l$ the $L[0]$ eigen-value of $s[0]^lv$ is too
small to fit into the subspace based on $|m_{\rm original}+lm,n\vac$,
see Remark \ref{lflat}.

Since $\VOA_{C,g}$ has structure of loop-module over $\CC[C]$,
its localization has structure of loop-module over $\CC[C_1]$.
One can also show that vertex algebra structure on $\VOA_{C_1,g}$
is the localization of the structure on $\VOA_{C,g}$.
\hfill $\Box$

\begin{rem}
{\rm 
It is interesting to observe that a surjective map on 
Fock spaces leads to an injective map on BRST cohomology.
}
\end{rem}

\begin{rem}
{\rm Even though $\ee^{\int m\cdot B(z)}$ is invertible,
its zero mode is not. This seems to contradict the calculations
of Section 2 but the reason is the presence of negative 
$L[0]$ eigen-values in $\Fock_{M\oplus K^*}$.
}
\end{rem}

So we now have in our disposal a way of calculating sections of \MSV\ of
a smooth toric variety on any
toric affine subset of it. This allows us to calculate
cohomology of chiral de Rham complex. 
We will be most concerned with a calculation of cohomology 
of \MSV\ for canonical line bundle $L$ over a complete toric variety
$\PP$. To get the fan of $L$ from the fan of $\PP$ one adds extra
dimension to $N$ and then lifts the fan of $\PP$ to height one as
illustrated  by the following figure.
$$
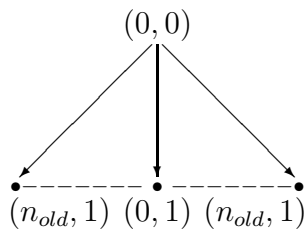
\begin{array}{ccc}
&\hspace{-50pt}\stackrel{~~~}{(0,0)}&\\[-10pt]
&
{\vector(-1,-1){50}\hspace{42pt}\vector(0,-1){50}\hspace{2pt}
\vector(1,-1){50}}&\\
\hspace{-50pt}\stackrel{\bullet-----}{(n_{old},1)}
&\hspace{-50pt}\stackrel{--\hspace{3pt}\bullet\hspace{3pt}--}{(0,1)}&
\hspace{-50pt}\stackrel{-----\bullet}{(n_{old},1)}\\
\end{array}
$$
%{\bf this was Figure 1.}

%calculation for the toric bundle
We adjust our notations to denote the whole new lattice 
by $N=N_1\oplus \ZZ$ and the new fan by $\Sigma$. An element
in $M$ that defines the last coordinate in $N$ is denoted
by "${\rm deg}$". Notice that for every cone $C^*\in \Sigma$
it is the same as "deg" from Proposition \ref{dimany}.
We also denote by $K^*$ the union of all cones in $\Sigma$,
which may or may not be convex.

As it was noticed in Remark \ref{pushforward}, we can consider
quasi-loco sheaf $\pi_* \msv(L)$ on $\PP$. 
By Proposition \ref{acyclic}, its cohomology could be calculated as
\v{C}ech cohomology that corresponds to the covering of
$L$ by open affine subsets $\AA_C$ where we only consider the
cones $C^*$ that contain $(0,1)$. So we need to consider
\v{C}ech complex
$$0\to \oplus_{C^*_0} \VOA_{C_0,g}
\to  \oplus_{(C^*_0,C^*_1)} \VOA_{C_{01},g} \to...
\to \oplus_{(C^*_0,...,C^*_r)} \VOA_{C_{0...r},g} \to 0$$
where $C_{0..k}$ is the dual of the intersection $C^*_{0...k}$ of
$C_0^*,...,C_k^*$. Here we have chosen
non-zero numbers $g_n$ for all generators $n$ of one-dimensional
cones in $\Sigma$. We know that each $\VOA_{C,g}$ is BRST cohomology of
the corresponding Fock space and our goal is to write cohomology
of \v{C}ech complex as certain BRST cohomology.

\begin{prop}
{\rm 
Consider the following double complex
$$
\begin{array}{ccc}
 0          & ...              & 0 \\
 \downarrow & ...              & \downarrow\\
0\to\oplus_{C^*_0} (\Fock_{M\oplus C^*_0})_{{\rm deg}\cdot
B[0]=0}
 &\to ... \to&
\oplus_{(C^*_0,...,C^*_r)} (\Fock_{M\oplus C^*_{0...r}})_{{\rm deg}\cdot
B[0]=0}\to 0 \\
 \downarrow & ...              & \downarrow\\
0\to\oplus_{C^*_0} (\Fock_{M\oplus C^*_0})_{{\rm deg}\cdot
B[0]=1}
 &\to ... \to&
\oplus_{(C^*_0,...,C^*_r)} (\Fock_{M\oplus C^*_{0...r}})_{{\rm deg}\cdot
B[0]=1}\to 0 \\
 \downarrow & ...              & \downarrow\\
  ...& ...              & ...\\
\end{array}
$$
where vertical arrows are ${\cal BRST}_g$ operators and horizontal 
arrows are sums of surjective maps of Fock spaces as dictated 
by definition of \v{C}ech cohomology. We also multiply vertical
differentials in odd-numbered columns by $(-1)$ to assure 
anticommutation of small squares. Then $p$-th cohomology of total 
complex is equal to $H^p(\pi_*\msv(L),\PP)$. Here again we only 
consider cones $C^*$ that contain $(0,\RR_{\geq 0})$.
}
\label{double}
\end{prop}

{\em Proof.} Proposition \ref{dimany} tells us that cohomology 
along vertical lines happens only at the top (zeroth) row,
where it becomes the \v{C}ech complex for the 
sheaf $\pi_*\msv(L)$. So spectral sequence of one stupid 
filtration degenerates and converges to cohomology of
$\pi_*\msv(L)$.
\hfill${\Box}$

Our next step is to calculate cohomology of total complex 
using the other stupid filtration. Let us see what happens
if we take cohomology of horizontal maps of our
double complex first. It could be done separately for each 
lattice element $n\in N$. If ${\rm deg}\cdot n=l$ then we are
be dealing with the $l$-th row. The part of the complex that 
we care about is a constant space $\Fock_{M\oplus n}$
multiplied by a certain finite complex of vector spaces.
That complex calculates cohomology of the simplex based on all
indices $i$ such that cones $C^*_i$ that contain $(0,\RR_{\geq 0})$ and
$n$. If $n\notin K^*$ then the set is empty, and cohomology is zero.
However, if $n\in K^*$ the cohomology is $\CC$ and is located
at the zeroth column. As a result, horizontal cohomology is zero except 
the zeroth column. Therefore, cohomology could be calculated
by means of the restriction of ${\cal BRST}_g$ operator applied
to kernels of horizontal maps from the zeroth column.
The following theorem describes this space.

\begin{theo}
{\rm
Consider the following degeneration of the vertex algebra structure
on $\Fock_{M\oplus K^*}$. In the definition of vertex operator 
$\ee^{\int (m\cdot B+n\cdot A)(z)}$ when applied to $...|m_1,n_1\vac$,
the result is put to be zero, unless there is a cone in $\Sigma$ that 
contains both $n_1$ and $n$. This does provide a consistent set of
operator product expansions and the new algebra is denoted by
$\Fock_{M\oplus K}^\Sigma$.
We denote by $\Delta^*$ the set of all generators of one-dimensional
cones of $\Sigma$. We construct a BRST operator on $\Fock_{M\oplus K}^\Sigma$
by the formula 
$${\cal BRST}_g=\oint {\cal BRST}_g(z)dz = 
\oint \sum_{n\in \Delta^*} g_n (n\cdot\Psi)(z)\ee^{\int n\cdot A(z)}.$$
Then we claim that $\oplus_p H^p(\pi_*\msv(L))$ equals BRST cohomology
of $\Fock_{M\oplus K^*}^\Sigma$ with respect to ${\cal BRST}_g$.
}
\label{toricbundle}
\end{theo}

{\em Proof.} In view of Proposition \ref{double}, it is enough
to show that horizontal cohomology of the double complex of 
\ref{double} at zeroth column and the corresponding vertical differential
coincide with $\Fock_{M\oplus K^*}^\Sigma$ and ${\cal BRST}_g$.
Kernel of horizontal map consists of collections of elements
of $\Fock_{M\oplus C^*}$ that agree with restrictions. This can certainly
be identified with $\Fock_{M\oplus K^*}^\Sigma$ as follows.
For every point $n\in N$ we take the corresponding $n$-part of the
above collection of elements, since it is the same no mater which 
$C^*\ni n$ we choose. In the opposite direction, for each cone $C^*$ 
we take a sum of $n$-parts for all $n$ that belong to $C^*$. When
we apply vertical arrows to such collections of elements, for each $C^*$
we use only the part of ${\cal BRST}_g$ that contains $B[0]$ eigen-values
from that $C^*$. Under our identification this is precisely the
action of
the whole ${\cal BRST}_g$ on $\Fock_{M\oplus K^*}^\Sigma$ because
as a result of that action for any $n\in C^*$ the only terms that
survive and have a non-trivial projection back to $C^*$ come
from applying the part of ${\cal BRST}_g$ with $n$ in $C^*$.
\hfill{$\Box$}

We also want to show that the structure of vertex algebra induced
on BRST cohomology of $\Fock_{M\oplus K^*}^\Sigma$ coincides
with the vertex algebra structure on cohomology of 
$\pi_*\msv(L)$ defined in Proposition \ref{vertcoh}.
\begin{prop}
{\rm
Two structures of vertex algebra on $H^*(\pi_*\msv(L))$
coincide.
}
\label{samevert}
\end{prop}

{\em Proof.} The cup-product $(n)$ is induced on \v{C}ech 
cohomology by the following product on \v{C}ech cochains.
To define \v{C}ech differential we have chosen an order on 
the set of all cones. If 
$\alpha\in \VOA_{C_{0...k},g}$ 
where $C_0<C_1<...<C_k$, and 
$\beta\in \VOA_{C'_{0...l},g}$ 
where $C'_0<C'_1<...<C'_l$, then their $(n)$-product 
$\alpha_{(n)}\beta$ is zero unless 
$$
C_0<...<C_k=C'_0<...<C'_l,
$$
in which case it is defined as the $(n)$-product of the 
restrictions of $\alpha$ and $\beta$ to $\VOA_{C_{0...k}C'_{1...l},g}$.
We extend this construction to define
$\alpha_{(n)}\beta$ for any pair of elements of the double complex of Theorem
\ref{toricbundle} by replacing the $(n)$-product in $\VOA_{C_{0...k}C'_{1...l},g}$
by $(n)$-product in $\Fock_{M\oplus (C^*_{0...k}\cap C'^*_{1...l})}$.
We now observe that for the differential $d$ of the total
complex we have
$$d(\alpha_{(n)}\beta)=(d\alpha)_{(n)}\beta + (-1)^{{\rm parity}(\alpha)+{\rm
column}(\alpha)}
\alpha_{(n)}(d\beta).$$
To check this we again use equation 4.6.9 of 
\cite{Kac}. The product $(n)$ induces the product
on cohomology of vertical maps that is precisely the
$(n)$ product on \v{C}ech cochains of
$\pi_*\msv(L)$. It also induces a cup product on
the cohomology of horizontal maps. A map between
the two repeated cohomologies could be seen on the level 
of cochains as an addition of a coboundary, which is
therefore compatible with $(n)$. It remains to notice
that the $(n)$ product on the zeroth column of our double
complex simply acts as an independent application of $(n)$
products for every cone $C^*\in \Sigma$. So it coincides on
the cohomology of the horizontal maps with the $(n)$ product
of the vertex algebra structure of $\Fock_{M\oplus K^*}^\Sigma$.
\hfill$\Box$

\begin{rem}
{\rm 
It is important to keep in mind that operations ${(n)}$ do not
define the structure of vertex algebra on the whole double
complex, they only induce this structure on cohomology.
This is analogous to the fact that the usual cup-product is
not super-commutative on the level of cochains.
Also, we can not really define a quasi-loco sheaf
of vertex algebras  over $\PP$ whose sections over $\AA_C$ 
are $\Fock_{M\oplus C^*}$, because eigen-values of $L[0]$ 
are not bounded from below. Perhaps, it is just a matter 
of definitions, but localization might indeed behave poorly
in this case.
}
\end{rem}

\begin{rem}
{\rm 
It is clear that fields $L(z)$, $J(z)$, $G(z)$ and $Q(z)$
are still given by
$$
Q(z)=A(z)\cdot\Phi(z)-{\rm deg}\cdot\partial_z
\Phi(z),~G(z)=  
B(z)\cdot\Psi(z),
$$
$$
J(z)=\,\no\Phi(z)\cdot\Psi(z)\no + \,{\rm deg}\cdot B(z)
$$
$$L(z)=\,\no B(z)\cdot A(z)\no +\no \partial_z\Phi(z)\cdot
\Psi(z)\no.$$
} 
\label{ljgqtoric}
\end{rem}

\begin{rem}
{\rm 
Similar results can be obtained for cohomology of
\MSV\ for an arbitrary smooth toric variety, for example
for a projective space. However in this paper we are mostly concerned
with line bundle case because of its applications to Mirror 
Symmetry.
}
\end{rem}

So now we have a good description of cohomology of \MSV\ on
a canonical bundle of a smooth toric variety. Our next step
is to use the results of Section 5 to obtain a similar
result for a Calabi-Yau hypersurface inside a smooth toric
nef-Fano variety.

Certain combinatorial conditions on $\Sigma$ are necessary to
ensure  that we have a Calabi-Yau hypersurface. 
Details could be found in the original paper of Batyrev
\cite{bat.dual}. The set $\Delta^*$ of
all $n$ should be the set of all lattice points of a convex
polytope which we also denote by $\Delta^*$ abusing notations
slightly. We do not require that all $n$ except $(0,1)$ are vertices 
of $\Delta^*$ which geometrically means that the opposite of the canonical
divisor on $\PP$ is nef but not necessarily ample. We also have
a polytope $\Delta\in M$ defined as follows. Decomposition
$N=N_1\oplus \ZZ$ implies $M=M_1\oplus \ZZ$. We define
$$\Delta=K\cap \{M_1,1\}.$$
All vertices of polytope $\Delta$ belong to $M$. This is not
entirely obvious, but follows from the fact that all cones
of $\Sigma$ are basic and therefore $\Delta$ is reflexive in $M_1$.
A Calabi-Yau hypersurface in $\PP$ is given by a section of the negative
canonical line bundle of $\PP$. Any such section is given by a set 
of numbers $f_m$, one for each lattice point $m$ in $\Delta$.
If $f$ is generic, the resulting hypersurface in $\PP$ is smooth
and Calabi-Yau. In what follows we will denote $(0,1)$ by
${\rm deg}^*$.

Let us  fix a generic section $\mu$ of the anti-canonical line bundle
of $\PP$ and
hence a function $f:\Delta\to\CC$.
\begin{prop}
{\rm
Operator ${\cal BRST}_\mu(z)$ from Lemma \ref{dmu} is
given by 
$${\cal BRST}_\mu(z)=\sum_{m\in\Delta} f_m (m\cdot \Phi)(z)
\ee^{\int m\cdot B(z)}.$$
}
\label{bf}
\end{prop}

{\em Proof.} It is enough to consider a section $\mu=f_m x^m$.
For any maximum cone of $\Sigma$ with basis $(n_1,...,n_{\dim N})$ we
see that 
$${\cal BRST}_{\mu}(z) = f_m \sum_i \varphi_i 
(m\cdot n_i)\ee^{\int
(m-m_i)\cdot B(z)} 
=f_m \sum_i  (m_i\cdot \Phi)(z) (m\cdot n_i) \ee^{\int m\cdot B(z)}
$$
$$=
f_m (m\cdot \Phi)(z) \ee^{\int m\cdot B(z)}.$$
\hfill$\Box$

%calculation of the affine piece of the hypersurface
From now on we denote ${\cal BRST}_\mu$ by ${\cal BRST}_f$.
We also denote by $X$ the Calabi-Yau hypersurface in $\PP$
which is given by $f$. We will denote by $X_C$ the
intersection of $X$ with $\AA_C$.
\begin{prop}
{\rm
For every cone $C^*\in\Sigma$ sections of $\msv(X)$ on $X_C$ are
given by BRST cohomology of $\Fock_{M\oplus C^*}$ by the 
operator 
$$\Bfg=\oint \Bfg(z) dz$$
where 
$$\Bfg(z)={\cal BRST}_f(z)+{\cal BRST}_g(z)
$$
$$=
\sum_{m\in\Delta} f_m(m\cdot \Phi)(z) \ee^{\int m\cdot B(z)}+
\sum_{n\in\Delta^*\cap C^*} g_n(n\cdot \Psi)(z) \ee^{\int n\cdot A(z)}.$$
}
\label{affbfg}
\end{prop}

{\em Proof.} One easily computes that all modes of ${\cal BRST}_f(z)$
and ${\cal BRST}_g(z)$ anti-commute with each other. Also, Proposition
\ref{affine} implies that sections of $\msv(X)$ are
cohomology with respect to ${\cal BRST}_f$
of cohomology of $\Fock_{M\oplus C^*}$ with respect to 
${\cal BRST}_g$. Consider the following double complex. 
$$
\begin{array}{ccccccc}
...~~~&0&&0&&0&~~~...\\ 
      &\downarrow&     &\downarrow&    &\downarrow&   \\
...\to&\Fock_{-1,0}&\to&\Fock_{0,0}&\to&\Fock_{1,0}&\to...\\ 
      &\downarrow&     &\downarrow&    &\downarrow&   \\
...\to&\Fock_{-1,1}&\to&\Fock_{0,1}&\to&\Fock_{1,1}&\to...\\ 
      &\downarrow&     &\downarrow&    &\downarrow&   \\
...\to&\Fock_{-1,2}&\to&\Fock_{0,2}&\to&\Fock_{1,2}&\to...\\ 
      &\downarrow&     &\downarrow&    &\downarrow&   \\
&...&&...&&...&\\
\end{array}
$$
where $\Fock_{k,l}$ is a shorthand for the part of $\Fock_{M\oplus C^*}$
where $({\rm deg}^*\cdot A)[0]$ and $({\rm deg}\cdot B)[0]$ equal
$k$ and $l$ respectively. Horizontal maps are ${\cal BRST}_f$
and vertical maps are ${\cal BRST}_g$. We already know that
columns of this double complex are exact everywhere except 
the zeroth row. A standard diagram chase then implies that 
horizontal cohomology of zeroth kernels of vertical maps
are isomorphic to cohomology of total complex. However,
total complex and differential on it are precisely 
$\Fock_{M\oplus C^*}$  and $\Bfg$. 
\hfill$\Box$

\begin{prop}
{\rm 
In the above proposition all cohomology of total complex 
are trivial, except for the zeroth one.
}
\label{zerototal}
\end{prop}

{\em Proof.} Grading by ${\rm deg}^*\cdot A[0]$ on 
$\pi_*\msv(L)$ corresponds to counting $\sharp (b^1)+\sharp(\varphi^1)-
\sharp (a_1) - \sharp \psi_1$ where $x^1$ is the special coordinate
of the line bundle. Since this count is zero for $\msv(X)$, the
cohomology of ${\cal BRST}_f$ is concentrated at zeroth column.
\hfill$\Box$

\begin{rem}
{\rm
Above identification is also compatible with vertex algebra
structures. Really, this structure is induced from that of $\Fock_{M\oplus
C^*}$ both for the repeated and for the single use of BRST cohomology.
}
\end{rem}

%calculation for the hypersurface
Now we are in position to calculate cohomology of \MSV\
of Calabi-Yau hypersurfaces in toric Fano varieties.

\begin{theo}
{\rm
BRST cohomology
of $\Fock_{M\oplus K^*}^{\Sigma}$ with respect to 
BRST operator $\Bfg$ equals $H^*(X,\msv(X))$.
}
\label{conebig}
\end{theo}

{\em Proof.} Argument is completely analogous to that
of Proposition \ref{double}.
We construct a double complex similar to that of \ref{double},
but with ${\cal BRST}_g$ changed to ${\cal BRST}_{f,g}$
and ${\rm deg}\cdot B[0]$ changed to ${\rm deg}\cdot B[0]+{\rm deg}^* 
\cdot A[0]$. Proposition \ref{zerototal} assures that spectral
sequence of this double complex degenerates. \hfill$\Box$

It is now a technical matter to calculate fields $L(z)$,
$J(z)$, $G(z)$ and $Q(z)$.
\begin{prop}
{\rm
Fields $L$, $J$, $G$, and $Q$ on $H^*(\msv (X))$
are induced from the following fields on $\Fock_{M\oplus K^*}^{\Sigma}$
$$
Q(z)=A(z)\cdot\Phi(z)-{\rm deg}\cdot\partial_z \Phi(z)
$$
$$
G(z)= B(z)\cdot\Psi(z)-{\rm deg}^*\cdot\partial_z \Psi(z)
$$
$$
J(z)=\,\no\Phi(z)\cdot\Psi(z)\no +\, {\rm deg}\cdot B(z)
-{\rm deg}^*\cdot A(z)
$$
$$
L(z)=\,\no B(z)\cdot A(z)\no +\no \partial_z\Phi(z)\cdot\Psi(z)\no
-{\rm deg}^*\cdot\partial_z A(z).
$$
}
\label{jlgq}
\end{prop}

{\em Proof.} By a standard application of Wick theorem, one
observes that the above operators satisfy OPEs of the
topological algebra of dimension $(\dim N-2)$. Then it remains to 
show that $G(z)$ and $Q(z)$ are correct. This
follows from Proposition \ref{GZhyper} and Remark \ref{ljgqtoric}.
\hfill{$\Box$}

\begin{rem}
{\rm
Notice that OPEs of topological algebra hold exactly in 
$\Fock_{M\oplus K^*}^\Sigma$ even though we only need them
to hold modulo the image of $\Bfg$. This algebra was discovered
almost three years ago as a lucky guess motivated
by Mirror Symmetry. The above algebra behaves exactly like
the holomorphic part of $N=(2,2)$ theory under mirror
involution, see for example \cite{witten}. This
becomes apparent when we undo the topological twist and consider $N=2$
super-conformal algebra with 
$$
L_{N=2}(z)=\,\no B(z)\cdot A(z)\no +\no
(1/2)(\partial_z\Phi(z)\cdot\Psi(z)-\Phi\cdot\partial_z\Psi(z))\no
$$
$$
-(1/2){\rm deg}^*\cdot\partial_z A(z)-(1/2){\rm deg}\cdot\partial B(z).
$$
}
\end{rem}

Mirror Symmetry in this setup means switching $M$ and $N$, $\Delta$
and $\Delta^*$, ${\rm deg}$ and ${\rm deg}^*$, $f$ and $g$.
We want to show that A model vertex algebra of a Calabi-Yau
hypersurface equals B model algebra of its mirror.
Certainly the formulas above show that this is true for the
operators $G$, $Q$, $J$ and $L$. However,
there are some obstacles that prevent us from making such a statement. 
The easiest objection is that the cocycle used to make vertex operators
bosonic does not look symmetric. However, that could be fixed by noticing
that multiplication of $|m,n\vac$ by $(-1)^{m\cdot n}$ is equivalent
to switching the roles of $M$ and $N$ in the definition of the cocycle.
We also notice that when we switch $M$ and $N$, and
consequently $\Psi$ and $\Phi$, the moding of the operators changes
slightly. Besides, the Fock space we consider is based on $M\oplus
K^*$, rather than on $M\oplus N$ or $K\oplus K^*$. Both these objections
will be addressed successfully in the next section. However, there
is one more difficulty that can not be resolved: we use the subdivision
of the cone $K^*$ but we do not subdivide lattice $M$ at all. This 
will eventually lead to the following statement\\[5pt]
{\it A  and B model vertex algebras of mirror symmetric Calabi-Yau
hypersurfaces are different degenerations of two families of 
topological vertex algebras that are related by a mirror involution.}

\section{Transition to the whole lattice}
The goal of this section is to show that cone $K^*$ in Theorem
\ref{conebig} could be replaced by the whole lattice $N$. This is
certainly far from obvious. We construct and use certain
homotopy operators whose anticommutators with $\Bfg$ are identity
plus operators that "push elements closer to $K^*$".
Before we start, notice that decomposition of $K^*$ into 
cones can be extended to decomposition of lattice $N$
by adding arbitrary multiples of ${\rm deg}^*$ to all cones.
This allows us to define vertex algebra $\Fock_{M\oplus N}^\Sigma$
analogously to the definition of the vertex algebra 
$\Fock_{M\oplus K^*}^\Sigma$.

\begin{prop}
{\rm For every vertex $m_0$ of $\Delta$ there is an operator 
${\cal R}_{m_0}$ such that 
$${\cal R}_{m_0}\Bfg+\Bfg {\cal R}_{m_0}={\bf 1}+\alpha$$
where $\alpha$ strictly increases eigen-values of $m_0\cdot
B[0]$ and does not decrease eigen-values of $m\cdot B[0]$ for any other
$m\in\Delta$.
}
\label{homop}
\end{prop}
{\em Proof.}
Consider graded ring $\CC[K]$. Pick a basis $n_1,...,n_{\dim N}$
of $N$. It was proved in \cite{locstring} that for general values of
$f_m$ elements $\sum_m f_mx^m (m\cdot n_i) $ form a regular sequence.
In particular, quotient ring is Artinian and  for a sufficiently 
big $k$ ($k=\dim M$ is in fact always enough) element $x^{km_0}$
lies in the ideal generated by the above regular sequence.
So we have 
$$x^{km_0}=\sum_i h_i(x^\Delta)\sum_m f_mx^m(m\cdot n_i).$$
We now consider the following field
$${\cal R}_{m_0}(z)=\ee^{-\int km_0\cdot B(z)}
\sum_ih_i(\ee^{\int\Delta\cdot B})n_i\cdot \Psi(z).$$
We have the operator product expansion
$${\cal R}_{m_0}(z){\cal BRST}_f(w)
\sim
{1\over {z-w}}\ee^{-\int km_0\cdot B}\sum_ih_i(\ee^{\int \Delta\cdot B})
\sum_m f_m \ee^{\int m\cdot B}(m\cdot n_i)
={1\over{z-w}}.$$

With the usual abuse of notation, we introduce 
$${\cal R}_{m_0}=\oint{\cal R}_{m_0}(z)dz.$$
We  now argue that this operator satisfies the claim of this
proposition.
First of all, the above OPE shows us that 
$${\cal R}_{m_0}{\cal BRST}_f+{\cal BRST}_f{\cal R}_{m_0}={\bf 1}.$$
Let us look at its anticommutator with ${\cal BRST}_g$.
Operator product expansion of ${\cal R}_{m_0}(z)\ee^{\int n\cdot
A(w)}n\cdot \Psi(w)$ is non-singular if $m_0\cdot n=0$. Otherwise, 
fields in the OPE shift eigen-values of $m_0\cdot B[0]$ positively
(more precisely by $m_0\cdot n$). Of course, for any other $m\in \Delta$
eigen-values of $m\cdot B[0]$ are shifted by $m\cdot n$ which is
non-negative.
\hfill$\Box$

The above proposition provides us with necessary tools to prove the 
main result of this section. 
\begin{prop}
{\rm 
Cohomology of ${\Fock}_{M\oplus K^*}^\Sigma$ 
with respect to $\Bfg$ is isomorphic to cohomology of 
${\Fock}_{M\oplus N}^\Sigma$ with respect to $\Bfg$.
}
\label{wholeN}
\end{prop}

{\em Proof.} It is clear that $\Bfg {\Fock}_{M\oplus K^*}^\Sigma
\subseteq {\Fock}_{M\oplus K^*}^\Sigma$. Then one needs to prove the 
following two inclusions \\
\noindent
$\bullet$ ${\rm Ker}({\Fock}_{M\oplus K^*}^\Sigma)
\cap {\rm Im}({\Fock}_{M\oplus N}^\Sigma)\subseteq {\rm Im}
({\Fock}_{M\oplus K^*}^\Sigma);$\\
\noindent
$\bullet$ ${\rm Ker}({\Fock}_{M\oplus K^*}^\Sigma)
+ {\rm Im}({\Fock}_{M\oplus N}^\Sigma)\supseteq {\rm Ker}
({\Fock}_{M\oplus K^*}^\Sigma).$\\
\noindent
Notice that the opposite inclusions are obvious.

{\em First inclusion.} Assume that there exists an element
$v\in \Fock_{M\oplus K^*}^\Sigma$ such that $v=\Bfg v_1$ where $v_1$
does not lie in $\Fock_{M\oplus K^*}^\Sigma$. Moreover, of all such $v_1$
we pick the one which is "the closest to $\Fock_{M\oplus
K^*}^\Sigma$".
The distance is defined as follows. We look at all codimension
one faces of $K^*$ or equivalently all vertices of $\Delta$.
For every vertex $m$ of $\Delta$ we look at the maximum
eigen-value of $-m\cdot B[0]$ on components of $v_1$. We call the 
maximum of this number and zero the $m$-distance from $v_1$ to 
$\Fock_{M\oplus K^*}^\Sigma$. Then the total distance from $v_1$ to 
$\Fock_{M\oplus K^*}^\Sigma$ is the sum of $m$-distances for all vertices
$m$
of $\Delta$. So we pick $v_1$ with a minimum distance
and our goal is to show that this distance is zero.

If the distance is not zero, then for one of the vertices $m$ 
there is a component of $v_1$ with a negative eigen-value of $m\cdot
B[0]$.
We now apply the result of Proposition \ref{homop}. Consider
operator ${\cal R}_m$. We have
$$({\cal R}_m\Bfg+\Bfg{\cal R}_m)v=v+\alpha v.$$
So 
$$v=\Bfg({\cal R}_mv+\alpha v_1),$$
because $\alpha$ commutes with $\Bfg$. Notice now that 
${\cal R}_mv\in \Fock_{M\oplus K^*}^\Sigma$ and the distance from
$\alpha v_1$ to $\Fock_{M\oplus K^*}^\Sigma$ is strictly less than
the distance from $v_1$ to it. Really, the $m$-distance is smaller,
and $m_1$-distance is not bigger for any other vertex $m_1$ of $\Delta$.
This contradicts minimality of $v_1$.

{\em Second inclusion.} Our argument here is similar. Let $v$ be an
element
of $\Fock_{M\oplus N}^\Sigma$ such that $\Bfg v=0$. Then for every vertex
$m\in \Delta$ we have 
$$\Bfg{\cal R}_m v=v + \alpha_mv$$
and therefore
$$v\equiv -\alpha_mv\hspace{-5pt}
\pmod {{\rm Im}(\Fock_{M\oplus N}^\Sigma)}.$$
By applying $\alpha_m$ for different $m$ sufficiently many times,
we can again push $v$ into $\Fock_{M\oplus K^*}^\Sigma$.
\hfill$\Box$

We now combine Propositions \ref{wholeN}
and \ref{jlgq} with Theorem \ref{conebig} to formulate one of the main
results of this paper.
\begin{theo}
{\rm 
Let $X$ be a Calabi-Yau hypersurface in a smooth toric nef-Fano variety,
given by $f:\Delta\to\CC$ and a fan $\Sigma$. Then cohomology of
\MSV\ of $X$ equals BRST cohomology of $\Fock_{M\oplus N}^\Sigma$
by operator 
$$\Bfg=\oint(
\sum_{m\in\Delta} f_m(m\cdot \Phi)(z) \ee^{\int m\cdot B(z)}+
\sum_{n\in\Delta^*} g_n(n\cdot \Psi)(z) \ee^{\int n\cdot A(z)}
)dz$$
with any choice of non-zero numbers $g_n$.
Additional structure of topological vertex algebra is given by formulas
of Proposition \ref{jlgq}.
}
\label{maintheo}
\end{theo}

We have therefore addressed one of the questions posed at the end of last
section. Another obstacle for Mirror Symmetry stated there was the
fact that the modes of $\Phi$ and $\Psi$ are defined differently.
However, this is precisely what happens when we go from A model to B
model. Because $J=\no\Phi\cdot\Psi\no+{\rm deg}\cdot B-
{\rm \deg}^*\cdot A$, the moding of $\Psi$ and $\Phi$ changes when we
go from $L_{A-model}[0]$ to $L_{B-model}[0]=L_{A-model}[0]+J[0]$.
Really, while $\Phi[k]$ and $\Psi[k]$ change, the true
modes $\Phi_{(k)}$ and $\Psi_{(k)}$ are not affected by the switch
of the roles $M$ and $N$.

%degeneration
It remains to address the following question. What is the real 
meaning of going from $\Fock_{M\oplus N}$ to $\Fock_{M\oplus N}^\Sigma$?
It turns out that in the case when $\Sigma$ admits a convex
piece-wise linear function (which geometrically means that
$\PP$ is projective) this vertex algebra is a degeneration of
vertex algebra $\Fock_{M\oplus N}$.

The degeneration we are about do describe is completely analogous
to the one discussed in \cite{locstring} but is now performed for
the whole Fock space. Let $h:N_\RR \to \RR$ be a continuous function
which is linear on every cone of $\Sigma$ and satisfies
$$h(x+y)\leq h(x)+h(y)$$
with equality achieved if and only if $x$ and $y$ lie in the same
cone of $\Sigma$. Then we get ourselves a complex parameter $t$
and start changing the basis of $\Fock_{M\oplus N}$ by
assigning $|m,n\vac_t=t^{h(n)}|m,n\vac$. To preserve the definition
of vertex algebra we also multiply $\ee^{\int n\cdot A(z)}$
by $t^{h(n)}$. Now if we let $t$ go to zero, the structure of 
vertex algebra of $\Fock_{M\oplus N}$ will go to the structure
of vertex algebra of $\Fock_{M\oplus N}^\Sigma$. 
When the structure of $\Fock_{M\oplus N}^\Sigma$ is defined
via this limit, we can also get the action of $\Bfg$ on
it as a limit of 
$$
\Bfg(t)=\oint(
\sum_{m\in\Delta} f_m(m\cdot \Phi)(z) \ee^{\int m\cdot B(z)}+
\sum_{n\in\Delta^*} g_nt^{h(n)}(n\cdot \Psi)(z) \ee^{\int n\cdot
A(z)})dz
$$

This prompts the following definition.
\begin{dfn}
{\rm We define {\it Master Family of vertex algebras}\
that corresponds to the pair of reflexive polytopes $\Delta$
and $\Delta^*$ as the BRST quotient of
the vertex algebra $\Fock_{M\oplus N}$ by the operator
$$
\Bfg=\oint(
\sum_{m\in\Delta} f_m(m\cdot \Phi)(z) \ee^{\int m\cdot B(z)}+
\sum_{n\in\Delta^*} g_n(n\cdot \Psi)(z) \ee^{\int n\cdot
A(z)})dz
$$
where $f$ and $g$ are parameters of the theory. Additional
structure of  topological vertex algebra is given by formulas of
Proposition \ref{jlgq}.
}
\label{master}
\end{dfn}

\begin{conj}
{\rm 
Vertex algebras that appear in Mirror Symmetry for hypersurfaces
defined by $\Delta$ and $\Delta^*$ are degenerations of
Master Family of vertex algebras. 
}
\label{conjecture}
\end{conj}

\begin{rem}
{\rm
Large complex structure limit (see \cite{Morrison}) in our language is
most likely the double degeneration of Master Family where both $M$
and $N$ are subdivided. This theory is a degeneration of theories
on both mirror manifolds and can also be used to link the two.
}
\label{lcsl}
\end{rem}

Our discussion so far have been focused around reflexive polytopes
$\Delta^*$ that admit a unimodular triangulation and therefore
yield smooth $\PP$s. This is a very important class of examples,
which includes famous quintic in $\PP^4$, but
most reflexive polytopes do not fall into this category.
Next two sections will be devoted to the treatment of singular
$\PP$s. We can no longer use the definition of \cite{MSV}, but
many of our results still hold in that generality
under appropriate definitions.

\section{Vertex algebras of Gorenstein toric varieties}
The goal of this section is to define an analog of \MSV\ for
an arbitrary Gorenstein toric variety. It is again a loco sheaf of 
conformal vertex algebras. Sections of this  sheaf over any
toric affine chart admit a structure of topological vertex algebra which
may or may not be compatible with the localization. However, $J[0]$ and
$Q[0]$ are globally defined, which  allows us to introduce string de Rham
complex and to propose a definition of string cohomology vector
spaces. Recall that dimensions of these spaces were rigorously defined by
Batyrev and Dais in \cite{BD} but the spaces themselves have never been
constructed mathematically.

We are working in the following setup. 
There are dual lattices $M$ and $N$ with a primitive element
"${\rm deg}$" fixed in $M$. 
There is a fan $\Sigma$ in $N$ such that all generators
$n_i$ of its one-dimensional faces satisfy ${\rm deg}\cdot n_i=1$.
A set $\Delta^*$ consists
of some lattice points of degree one
inside the union of all cones of $\Sigma$. We do not generally
require that $\Delta^*$  includes all such points, or that 
it is a set of all lattice points inside a convex polytope. However,
we do demand that it contains generators of all one-dimensional 
cones of $\Sigma$. At last, we have a generic set of numbers $g_n$ for
all $n\in \Delta^*$.

\begin{dfn}
{\rm 
For each cone $C^*\in \Sigma$ we denote by ${\cal V}_g(C)$ the 
BRST cohomology of the vertex algebra $\Fock_{M\oplus C^*}$
with respect to the BRST operator 
$${\cal BRST}_g=\oint(\sum_{n\in \Delta^*\cap C^*} g_n
(n\cdot \Psi)(z)\ee^{\int n\cdot A(z)})dz.$$
We also provide this algebra with structure of topological
algebra by introducing operators
$$
Q(z)=A(z)\cdot\Phi(z)-{\rm deg}\cdot \partial_z\Phi(z),~G(z)=
B(z)\cdot\Psi(z),
$$
$$
J(z)=\,\no\Phi(z)\cdot\Psi(z)\no +\, {\rm deg}\cdot B(z),
$$
$$L(z)=\,\no B(z)\cdot A(z)\no +\no \partial_z\Phi(z)\cdot
\Psi(z)\no.$$
}
\label{localvoa}
\end{dfn}

\begin{rem}
{\rm 
Above definition does not guarantee that the resulting 
vertex algebra has no negative eigen-values of $L[0]$.
To prove this and much more we first consider the case
of simplicial cone $C^*$. Then we extend our results to the
general case by looking at the degeneration of $\Fock_{M\oplus C^*}$ that
corresponds to subdivision of $C^*$ into simplicial cones.
}
\end{rem}

First of all, we consider the case of orbifold singularities
in which case we can give an explicit description of 
BRST cohomology similar to that of Proposition \ref{dimany}.
When we talk about orbifold singularities we implicitly
assume that not only the cones $C^*$ are simplicial, but also
$g_n$ are zero, except for the generators of one-dimensional
faces of $C^*$.
\begin{prop}
{\rm 
Let $C^*$ be a simplicial cone of dimension $\dim N$. Its faces
of dimension one are generated by
$n_1,n_2,...,n_{\dim C^*}$. Denote by $N_{small}$ the sublattice of $N$
generated by $n_i$. Denote by $M_{big}$ the suplattice of $M$ which is
the dual of $N_{small}$. Let the dual of $C^*$ in $M_{big}$
be generated by $m_1,...,m_{\dim M}$. For every $i$ we define
$$b^i(z)=\ee^{\int m_i\cdot B(z)},
~\varphi^i(z)=(m_i\cdot \Phi(z))\ee^{\int
m_i\cdot B(z)}, ~\psi_i(z)=(n_i\cdot \Psi(z))\ee^{-\int m_i\cdot B(z)},$$
$$a_i(z)=\,\no (n_i\cdot A(z))\ee^{-\int m_i\cdot B(z)}\no +\no
(m_i\cdot\Phi(z))(n_i\cdot\Psi(z))\ee^{-\int
m_i\cdot B(z)}\no$$ 
for all $i=1,...,\dim M$. 
These fields generate a vertex subalgebra $\VOA_{C^*,M_{big}}$ inside
$\Fock_{M_{big}\oplus 0}$. Consider all fields from $\VOA_{C^*,M_{big}}$
whose $A[0]$ eigen-values lie in $M$. Denote the resulting algebra
by $\VOA_{C^*,M}$. 
Let ${\rm Box}(C^*)$ be the set of all elements
in $n\in C^*$ such that $(n-n_i)\notin C^*$ for all $i$.  
For every $n\in {\rm Box}(C^*)$ consider the following set of elements
of $\Fock_{M\oplus n}$. For every 
$v_0=\prod A[...]\prod B[...]\prod
\Phi[...]\prod\Psi[...]|m,0\vac$ 
that lies in $\VOA_{C^*,M}$
consider $v=\prod A[...]\prod B[...]\prod
\Phi[...]\prod\Psi[...]|m,n\vac$ which is obtained by applying the same
modes of  $A$, $B$, $\Phi$ and $\Psi$ to $|m,n\vac$ instead of $|m,0\vac$.
We denote this space by  
$\VOA_{C^*,M}^{(n)}$.
Then we claim that 
$${\cal V}_g(C) = \oplus_{n\in {\rm Box}(C^*)} \VOA_{C^*,M}^{(n)}.$$\\[-3em]
}
\label{orbiloc}
\end{prop}

{\em Proof.} First of all, the argument of Proportion \ref{dimone}
shows that 
$${\cal V}_g(C) = \oplus_{n\in {\rm Box}(C^*)} 
{\rm Ker}({\cal BRST}_g:\Fock_{M\oplus n}\to \Fock_{M\oplus C^*}).$$
Then we notice that changing $|m,n\vac$ to $|m,0\vac$ commutes with the 
action of ${\cal BRST}_g$ so it is enough to concern ourselves with
the case $n=0$. Then the kernel of ${\cal BRST}_g$ on 
$\Fock_{M\oplus 0}$ is the intersection of $\Fock_{M\oplus 0}$ 
with the kernel of ${\cal BRST}_g$ on $\Fock_{M_{big}\oplus 0}$.
It remains to apply Proposition \ref{dimany}.
\hfill${\Box}$

\begin{rem}
{\rm
Corresponding space $\AA_C$ is a quotient of a flat 
space by an abelian group. The part at $n=0$ is precisely
the invariant part of the flat space algebra, while other
$n$ correspond to "twisted sectors".
}
\end{rem}

\begin{prop}
{\rm 
For a simplicial cone $C^*$ of dimension $\dim N$ 
eigen-values of 
$L[0]$ on ${\cal V}_g(C)$ are non-negative.
Eigen-values of $A[0]$ on the zero eigen-space of $L[0]$
lie in $C$. Besides, for any $d>0$, eigenvalues of 
$A[0]$ on $L[0]=d$ eigen-space lie in $C-D(d)\,{\rm deg}$
where $D(d)$ is some constant which depends only on $d$ and
dimension of $N$.
}
\label{orbiL}
\end{prop}

{\em Proof.} Consider $n\in {\rm Box}(C^*)$
given by
$$n=\sum_i\alpha_i n_i.$$
Let us consider all elements $v_0$ from $\VOA_{C^*,M_{big}}$
and the corresponding elements $v$ obtained by changing
$|m,0\vac$ to $|m,n\vac$. Such a change incurs a change in
$L[0]$ which is equal to $m\cdot n$, and is therefore
linear in $m$. Hence, we can assess new $L[0]$
by adding $\alpha_i$ for each occurrence of $b^i$ or 
$\varphi^i$ and subtracting $\alpha_i$ for each occurrence
of $a_i$ or $\psi_i$. So to show that $L[0]$ has no negative
eigen-values, one has to show that each mode of $a_i$,
$b^i$, $\varphi^i$ or $\psi_i$ contributes non-negatively.
According to Proposition \ref{dimany}, we need to consider
non-positive modes of $b$ and $\varphi$ and negative modes 
of $a$ and $\psi$. 
The resulting contributions are collected in the following
table.
$$
\begin{array}{rcccc}
{\rm mode}~\vline &a_i[-k] &b^i[-k] &\varphi^i[-k] &\psi_i[-k]\\
\hline
\Delta L[0]~\vline&k-\alpha_i&k+\alpha_i&k+\alpha_i&k-\alpha_i\\
\end{array}
$$
Since $0\leq \alpha_i <1$, all entries are non-negative,
which proves the first part of the proposition.
Moreover, the only time a zero entry can occur is when
$\alpha_i=0$ and we are looking at $b^i[0]$ or $\varphi^i[0]$,
which proves the second part. Finally, for a fixed $d$ there
are only finitely many ways to combine $a[-k]$, $b[-k]$, $\varphi[-k]$ 
and $\psi[-k]$ to get $L[0]=d$, up to arbitrary extra
$b^i[0]$ and $\varphi^i[0]$. This finishes the proof.
$\hfill{\Box}$

We now extend these results to simplicial cones
of dimension smaller than $\dim N$.
\begin{prop}
{\rm 
For any simplicial cone $C^*$ 
eigen-values of 
$L[0]$ on ${\cal V}_g(C)$ are non-negative.
Eigen-values of $A[0]$ on the zero eigen-space of $L[0]$
lie in $C$. 
Besides, for any $d>0$, eigenvalues of 
$A[0]$ on $L[0]=d$ eigen-space lie in $C-D(d)\,{\rm deg}$
where $D(d)$ is some constant which depends only on $d$ and
dimension of $N$.
}
\label{anysimpl}
\end{prop}

{\em Proof.} Consider an arbitrary simplicial cone $C_1^*$
of maximum dimension whose one-dimensional faces are generated
by elements in $\Delta^*$ such that $C^*$ is a face of $C_1^*$.
Proposition \ref{orbiL} implies
that ${\cal V}_{C_1}$ is $\CC[C_1]$-loop-module. Proof 
of Proposition \ref{localization} is applicable in this more
general situation and allows us to show that ${\cal V}_C$
is a localization of ${\cal V}_{C_1}$ with respect to
the multiplicative system $S=\{x^m,m\cdot C^*=0,m\cdot C_1^*\geq 0\}$.
Grading operator $L[0]$ is still non-negative on the localization,
and $L[0]=0$ part of ${\cal V}_C$ is the localization of $L=0$
part of ${\cal V}_{C_1}$. It remains to observe that for every $x^m\in S$
we have  $-m\in C$ so localization does not push $A[0]$ eigen-values
from $C$. 
\hfill$\Box$

We are now in position to drop simpliciality assumption on  
cone $C^*$. This will require a careful investigation of degeneration
of vertex algebras given by a triangulation of a non-simplicial cone.

\begin{prop}
{\rm 
For any cone $C^*$ and a generic choice of $g$,
all eigen-values of $L[0]$ on ${\cal V}_C$ are non-negative.
}
\label{boundedL}
\end{prop}

{\em Proof.} Consider an arbitrary regular triangulation of 
$C^*\cap \Delta^*$ and corresponding decomposition of 
$C^*$ into union of simplicial cones.
Denote by $\Fock_{M\oplus C^*}^{triang}$ degeneration
of vertex algebra $\Fock_{M\oplus C^*}$ as in
Theorem \ref{toricbundle}. 
As in Proposition \ref{double} we consider the double complex
$$
\begin{array}{ccc}
 0          & ...              & 0 \\
 \downarrow & ...              & \downarrow\\
0\to\oplus_{C^*_0} (\Fock_{M\oplus C^*_0})_{{\rm deg}\cdot
B[0]=0}
 &\to ... \to&
\oplus_{(C^*_0,...,C^*_r)} (\Fock_{M\oplus C^*_{0...r}})_{{\rm deg}\cdot
B[0]=0}\to 0 \\
 \downarrow & ...              & \downarrow\\
0\to\oplus_{C^*_0} (\Fock_{M\oplus C^*_0})_{{\rm deg}\cdot
B[0]=1}
 &\to ... \to&
\oplus_{(C^*_0,...,C^*_r)} (\Fock_{M\oplus C^*_{0...r}})_{{\rm deg}\cdot
B[0]=1}\to 0 \\
 \downarrow & ...              & \downarrow\\
  ...& ...              & ...\\
\end{array}
$$
where $C_0^*$ are all cones in the triangulation. 
Again we consider two spectral sequences associated to this double
complex. When you take horizontal maps first, the only non-trivial
cohomology appears at zeroth column. Moreover,
when you apply vertical cohomology to zeroth kernels of horizontal maps,
you get BRST cohomology of $\Fock_{M\oplus C^*}^{triang}$
by ${\cal BRST}_g$. This shows that cohomology of the total
complex is isomorphic to BRST cohomology of $\Fock_{M\oplus
C^*}^{triang}$. The other stupid filtration implies that there
exists a spectral sequence that converges to BRST cohomology
of $\Fock_{M\oplus C^*}^{triang}$ and starts with
$\oplus_{C_{0...k}}{\cal V}_g(C_{0...k})$. Everything in this picture
is additionally graded by eigen-values of $L[0]$. Proposition
\ref{anysimpl} shows that 
${\cal V}_g(C_{0...k})$ have no negative eigen-values of $L[0]$,
therefore BRST cohomology of
$\Fock_{M\oplus C^*}^{triang}$ can not have any negative 
eigen-values either.

Now it remains to go from BRST cohomology of 
$\Fock_{M\oplus C^*}^{triang}$ to that of $\Fock_{M\oplus C^*}$.
Notice that $A[0]$, $L[0]$ and $J[0]$ commute with each other and 
with ${\cal BRST}_g$ so we can consider separately the parts of
$\Fock_{M\oplus C^*}$ and $\Fock_{M\oplus C^*}^{triang}$ that have fixed
eigenvalues of $A[0]$, $J[0]$ and
$L[0]$. We can show that these spaces are finite-dimensional.
Really, for a fixed eigenvalue $m$ of $A[0]$ we can find an integer $r$
such that $m+r\,{\rm deg}$ lies in $C$. Then we claim that all eigen-spaces
of $L[0]+(r+1)J[0]$ on
$\Fock_{m\oplus C^*}$ are finite-dimensional. For each $n$
we start with at least ${\rm deg}\cdot n$ and almost all modes of 
$A$, $B$, $\Phi$ and $\Psi$ can only increase this eigenvalue. The
exception is a few modes of $\Phi$ or $\Psi$, but they are 
fermionic and can only appear in a finite number of combinations.
This proves that all eigen-spaces of $L[0]+(r+1)J[0]$ are
finite-dimensional, and so are all spaces with fixed $A[0]$, $L[0]$
and $J[0]$ eigen-values.

Now ${\cal BRST}_g$ for $\Fock_{M\oplus C^*}^{triang}$
can be seen as the true limit of operators ${\cal BRST}_{g(t)}$
as discussed right before the Definition \ref{master}.
For families of operators on finite dimensional spaces, 
dimensions of kernels jump at special points and dimensions of
images decrease, so dimensions of cohomology jump.
Since for all negative $L[0]$ there is no cohomology for the
degenerate map, there is no cohomology for the original
$\Fock_{M\oplus C^*}$ for generic (that is outside of countably many
Zariski closed subsets of codimension one) choices of $g$.
\hfill$\Box$

We can use similar arguments to extend the rest of the results of 
Propositions \ref{orbiL} and \ref{anysimpl} to arbitrary Gorenstein
cones $C^*$.
\begin{prop}
{\rm
For any cone $C^*$,
the eigen-values of $A[0]$
on the zero eigen-space of $L[0]$  lie in $C$.
Besides, for any $d>0$, eigenvalues of 
$A[0]$ on $L[0]=d$ eigen-space lie in $C-D(d)\,{\rm deg}$
where $D(d)$ is some constant which depends only on $d$ and
dimension of $N$.
}
\label{lzero}
\end{prop}

{\em Proof.} For $L[0]=0$, it is enough to show for every generator
$n$ of a one-dimensional face of $C^*$ that eigen-values of $n\cdot A[0]$
are non-negative on the zero eigen-space of $L[0]$.
Argument of the above proposition shows that it
is enough to produce a regular triangulation of
$C^*$ such that all its maximum cones contain $n$. Really,
we can then apply the second part of Proposition
\ref{anysimpl} and the same degeneration argument works.
Similarly, to show that for $L[0]=d$ all eigenvalues
of $A[0]$ lie in $C-D(d)\,{\rm deg}$ means to
show that $n\cdot A[0]\geq -D(d)$ for every generator $n$ of a
one-dimensional face of $C^*$. Again if we can produce a regular
triangulation such that all of its maximum cones contain $n$, then 
we can use the same degeneration argument and Proposition  
\ref{anysimpl}.

To construct such a triangulation, we do the following.
We consider polytope $P=C^*\cap \Delta^*$. 
For every vertex $n_i$ of $P$ we move it slightly
away from $n$ along the line from $n$ to $n_i$. For
small generic perturbation of this type, the resulting
(non-integer) points $n_i$ and $n$ will still
be a set of vertices of a convex polytope $P'$.
All faces of $P'$ that do not contain $n$ will
be simplicial, and we will call the union of these
faces the outer surface of $P'$.
Then we consider the following function on
$h_1:P\to \RR$. For every point $p\in P$ we draw a line
$l$ from $n$ to $p$ and define $h(p)$ as the ratio of the
distances from $n$ to $p$ and from $n$ to the outer
boundary of $P'$. This function will be piecewise
linear on the triangulation of $P$ that is obtained 
by intersecting $P$ with ${\rm Conv}(n,F)$ for
all faces $F$ of the outer boundary of $P'$. Moreover,
it will be strictly convex in a sense that for
each two $p_1,p_2\in P$ and each $\alpha\in(0,1)$ 
$$\alpha h(p_1)+(1-\alpha)h(p_2)\leq h(\alpha p_1+
(1-\alpha)p_2)$$
with equality satisfied if and only if 
there exists a simplex of the triangulation that
contains both $p_1$ and $p_2$. We then extend $h$ from $P$ to 
$C^*$ by putting $h_{C^*}(p)=h(p)({\rm deg} \cdot  p)$.
This function will be strictly convex 
on the triangulation of $C^*$ such that all its simplices 
contain $n$, which finishes the proof.
\hfill$\Box$

\begin{rem}
{\rm
Looking back, we really had to work hard to prove 
last two propositions. It would be very interesting to
find a direct proof not based on results of
Proposition \ref{dimany}.
}
\end{rem}

\begin{prop}
{\rm
If $C_1^*$ is a face of $C^*$ then surjective map
$\Fock_{M\oplus C^*}\to \Fock_{M\oplus C_1^*}$
induces a map ${\cal V}_C\to {\cal V}_{C_1}$
which is a localization map of loop-module ${\cal V}_C$
over $\CC[C]$ by multiplicative system
$S=\{x^m,m\cdot C_1^*=0,m\cdot C^*\geq 0\}$.
}
\label{coneloc}
\end{prop}

{\em Proof.} The hard part was to show that 
spaces in question are loop-modules, that is to
show that $L[0]$ is non-negative. Then argument
of Proposition \ref{localization} extends to this
more general situation.
\hfill$\Box$

For the rest of this section we no longer assume that
${\rm deg}\in M$ is fixed. 

\begin{prop}
{\rm 
A toric variety $\PP$ defined by fan $\Sigma$
is Gorenstein if for every cone $C^*\in \Sigma$
all generators $n$ of one-dimensional cones of $C$ 
satisfy ${\rm deg}_C\cdot n=1$ where ${\rm deg}_C$
is a lattice point in $M$.
}
\end{prop}

{\em Proof.} See \cite{bat.dual} \hfill$\Box$

The following definition is made possible by the 
results of Propositions \ref{coneloc} and \ref{boundedL}.
\begin{dfn}
{\rm Let $\PP$ be a Gorenstein toric variety,
given by fan $\Sigma$ in $N$. Fix a generic set of 
numbers $g_n$ for all points of degree one 
in each cone of $\Sigma$ (the notion of degree
may vary from cone to cone). Then
(generalized) \MSV\ $\msv(\PP)$ is defined 
as a quasi-loco sheaf over it such that for any affine
subspace $\AA_C$ of $\PP$ sections of $\msv(\PP)$ are
BRST cohomology of $\Fock_{M\oplus C^*}$
by operator
$${\cal BRST}_g=\oint(\sum_{n\in C^*,{\rm deg}_C\cdot n=1}
g_n (n\cdot \Psi)(z)\ee^{\int n\cdot A(z)} )dz.$$
}
\label{defmsvgor}
\end{dfn}

\begin{rem}
{\rm 
Notice that while the choice of $g$ is irrelevant in
the smooth or even orbifold case, it is very important
in general.
}
\end{rem}

\begin{theo}
{\rm 
Quasi-loco sheaf
$\msv(\PP)$ is in fact a loco sheaf. Recall that the grading is given
by $L[0]$ where 
$$L(z)=\,\no B(z)\cdot A(z)\no +\no \partial_z\Phi(z)\cdot
\Psi(z)\no.$$
}
\label{proploco}
\end{theo}

{\em Proof.}
This is essentially a local statement, so it is enough to work
with one cone $C^*$ of maximum dimension.
We need a few preliminary lemmas.

\begin{lem}
{\rm 
The action of ${\rm deg} \cdot B [0]$ could be pushed down to the 
BRST cohomology of $\Fock_{M\oplus K^*}$. BRST cohomology of 
$\Fock_{M\oplus K^*}$ has only eigenvalues of ${\rm deg} \cdot B [0]$ 
in a certain range (from $0$ to $D_1$). 
}
\label{lemma2}
\end{lem}

{\em Proof of the Lemma.} Our BRST operator increases 
${\rm deg} \cdot B [0]$ eigen-values by one, so the action
of ${\rm deg} \cdot B [0]$ could be pushed down to the cohomology.
As before, we can notice that the above statement is true for 
any simplicial subcone and then do the spectral sequence and degeneration
trick as in Proposition \ref{boundedL}. We must mention that 
as a result of the spectral sequence the bound $D_1$ could jump a bit, 
but we only need to know that such a bound exists. $\hfill{\Box}$

We will need some general theory of funny objects which we call
{\it almost-modules}. 
\begin{dfn}
{\rm
Let $R$ be a Noetherian ring. An abelian group $V$ is called an {\it
almost-module}\ over $R$ if\\
\noindent$\bullet$
there is defined a map $R\times V\to V$ which is bilinear but not
necessarily associative;\\
\noindent$\bullet$
there exists a finite filtration of $V$ compatible with the multiplication
map above such that the quotients of the filtration are modules over $R$.
}
\end{dfn}

\begin{lem}
{\rm
Let $V$ be an almost-module over $R$. The following conditions are
equivalent:

(1) There exists a filtration of $V$ as above such that all quotients 
are finitely generated.

(2) For any filtration of $V$ as above all quotients 
are finitely generated.

(3) Any ascending chain of sub-almost-modules terminates.

If these conditions hold, then $V$ is called a Noetherian almost-module.
}
\label{nalmost}
\end{lem}

{\em Proof of the lemma.} 
$(3)\Rightarrow(2)$. If one of the quotients is not finitely generated,
then then there will be an ascending non-terminating chain of submodules
inside it. It will give rise to an ascending non-terminating chain of
sub-almost-modules of $V$.\\
$(1)\Rightarrow(3)$. If 
$$V=F^kV\subset F^{k-1}V \subset ...\subset F^0V=0$$ 
is the filtration on $V$ then for any ascending chain $\{V_j\}$ one can
look at $V_j\cap F^1V$. At some point it will stabilize. Then we can look
at $V_j\cup F^2V/F^1V$, and it will stabilize too. This will imply that 
$V_j\cup F^2V$ is stabilized. Eventually we will have $V_j\cup F^kV=V_j$
stabilized. $\hfill\Box$

As a corollary from the above proposition, all submodules and quotients of
Noetherian almost-modules are Noetherian.

\begin{lem}
{\rm
For any $d, D$ and any finite set $I\subseteq N$  the $L[0]=d$ part of
${\cal BRST}_g$ cohomology of $\Fock_{K-D\deg\oplus I}$ is a Noetherian
almost-module over $\CC[K]$. The multiplication is given by
$\ee^{\int m\cdot B}[0]$. 
}
\end{lem}

{\em Proof of the lemma.} Consider the filtration by the
number of $A$. When $\ee^{\int B\cdot m}[0]$ acts on the quotients
of this filtration, it amounts to acting on $|m_1,n_1\vac$ directly,
because one can push through and the extra commutators are in the lower 
part of the filtration. We can also do it for one fixed $n$, because there
are only finitely many of them.

Then we see that the multiplication gives zero unless $m\cdot n$=0,
and for those $m$ it is simply a shift of $A[0]$ eigen-values.
Now it remains to notice that there are finitely many (linearly
independent) choices of extra $A,B,\Phi,\Psi$ to get the desired $L[0]$ 
eigen-value. Also, there are finitely many choices for $m\cdot n$ in
$m,n\vac$ because it can;t be too big. And for each $m\cdot n=const$ 
we have a finitely generated module over $\CC[K]$ which finishes the
proof of the lemma. $\hfill\Box$

We are now ready to complete the proof of Theorem \ref{proploco}.
What we need to show is that $L[0]=d$ component of ${\cal BRST}_g$
cohomology of $\Fock_{M\oplus K^*}$ is a Noetherian almost module over
$\CC[K]$. Because of Proposition \ref{lzero} and Lemma \ref{lemma2}, it is
enough to consider $\Fock_{K-D\deg,\{\deg\cdot .\leq D_1\}}$.
By the above lemma, this is a Noetherian almost-module, and the
observation after Lemma \ref{nalmost} finishes the proof of the theorem.
$\hfill\Box$ 

The above theorem together with Proposition \ref{acyclic} lead to
the following corollary.
\begin{coro}
{\rm
If $\PP$ is compact, then $H^*(\msv(\PP))$ is a graded vertex
algebra with finite dimensional graded components.
}
\end{coro}

It is worthwhile to mention that $\msv(\PP)$ is also loop-coherent
with respect to the grading by the B model Virasoro operator $L_B[0]$
which is the zeroth mode of 
$$L_B(z)=\,\no B(z)\cdot A(z)\no +\no \partial_z\Psi(z)\cdot
\Phi(z)\no -{\rm deg}\cdot B(z).$$  
\begin{theo}
{\rm $\msv(\PP)$ is loop-coherent with respect to the grading by
$L_B[0]$.}
\label{Bproploco}
\end{theo}

{\em Proof.}
The proof is completely analogous to the proof of Theorem 
\ref{proploco}. We simply follow the chain of propositions 
of this section, with the following change. In the proof of
Proposition \ref{orbiL}, in addition to the table 
$$
\begin{array}{rcccc}
{\rm mode}~\vline &a_i[-k] &b^i[-k] &\varphi^i[-k] &\psi_i[-k]\\
\hline
\Delta L_B[0]~\vline&k-\alpha_i&k+\alpha_i&k+1+\alpha_i&k-1-\alpha_i\\
\end{array}
$$
we need to consider the extra shift by $\deg\cdot\alpha=\sum_i \alpha_i$
due to the extra term in $L_B(z)$. This shift cures possible negative
contributions of $\psi_i[-1]$ (fortunately, these are fermionic modes so
they can not repeat). $\hfill{\Box}$

We will now address the problem of 
{\it string-theoretic cohomology vector spaces}.
Recall that string cohomology numbers were constructed 
by Batyrev and Dais in \cite{BD}, and it was proved
in \cite{bat.bor} that they comply with predictions
of Mirror Symmetry for Calabi-Yau complete intersections
in toric varieties. Unfortunately, until now it was
not known how to construct string cohomology vector spaces
whose dimensions are the string cohomology numbers above.
The analysis of this paper suggests the following 
definition, at least for toric varieties.

\begin{dfn}
{\rm 
{\it String-differential forms}\ on $\PP$ is
$L[0]=0$ part of $\msv(\PP)$. By Remark \ref{zerocoh}
it is a coherent sheaf.
}
\end{dfn}

The following proposition 
provides us with a much more practical definition
of this sheaf that does not refer explicitly to sheaves
of vertex algebras.
\begin{prop}
{\rm 
For each $C^*\in \Sigma$ consider $\CC[C\oplus C^*]$-module
$V_C$ defined as 
$$V_C=\oplus_{m\in C,n\in C^*,m\cdot n=0}\CC\  x^my^n$$
where the action of $x^k$ and $y^l$ is defined to be zero if 
the result violates $m\cdot n=0$.
There is defined a differential 
${\cal BRST}_g^{(0)}$
on $V_C\otimes_{\CC}(\Lambda^*M_{\CC})$
given by 
$${\cal BRST}_g^{(0)} = \sum_{n\in C^*\cap \Delta^*}
g_n y^n {\rm contr}(n)$$
where ${\rm contr}(n)$ indicates contraction by $n$ on 
$\Lambda^*M$.
Then cohomology of ${\cal BRST}_g^{(0)}$ is
a finitely generated $\CC[C]$-module isomorphic
to $L[0]=0$ component of ${\cal V}_C$.
Moreover, grading by $J[0]$ on it is defined as 
"degree in $\Lambda^*M$ plus degree of $n$",
and differential $d$ is defined as
$$d(wx^my^n)=(w\wedge m)x^my^n.$$
}
\label{easyderham}
\end{prop}

{\em Proof.} Due to Proposition \ref{lzero}, $L[0]=0$ part 
of ${\cal V}_C$ could be obtained by applying ${\cal BRST}_g$
to $L[0]=0$ part of $\Fock_{C\oplus C^*}$. For every
$|m,n\vac$ from this space we already have $m\cdot n=0$,
so all elements from this space are obtained by multiplying 
$|m,n\vac$ by products of $\Phi^i[0]$, that is by $\Lambda^*M$.
Then we only need to calculate the action of ${\cal BRST}_g$
on this space, as well as the actions of $J[0]$ and $Q[0]$.
This is accomplished by a direct calculation.
\hfill$\Box$

If desired, one can use the above proposition as a definition
of the space of sections of the sheaf of string-differential
forms over $\PP$. It is not hard to show that it is coherent
directly.

\begin{rem}
{\rm
Notice that sheaf of string-differential forms is
not locally free, it reflects singularities of $\PP$. 
}
\end{rem}

\begin{rem}
{\rm
Another peculiar feature of the above description is that
grading by eigen-values of $J[0]$ on
the space of differential forms seems to be
ill-defined, since $J[0]$ varies with the cone.
Nevertheless, this is not a problem, because
the notion of $J[0]$ behaves well under the localization,
so string cohomology spaces do have an expected double grading.
Besides, de Rham operator $\oint Q(z)dz$ is clearly well-defined 
for string-differential forms on $\PP$.
}
\end{rem}

\begin{rem}
{\rm 
The fiber of the sheaf of string-differential forms over
the most singular point of $\AA_C$ is obtained 
by considering only $m=0$ part of the above space.
It is easily seen to coincide with the prediction of \cite{locstring}.
}
\end{rem}

\begin{rem}
{\rm 
It is not entirely clear if one should consider the cohomology
of $\msv(\PP)$ or the hypercohomology of it under the $Q(0)$
operator. On one hand, hypercohomology might be a smaller
and nicer object, but on the other hand taking hypercohomology
may complicate the relation between A and B models. So our definitions
below should be considered only provisional.
}
\end{rem}

\begin{dfn}
{\rm 
String cohomology vector space is the hypercohomology of the complex of
string-differential forms.
}
\end{dfn}

It is very likely that our definition reproduces
correctly the numbers of \cite{BD}, but clearly more
work is necessary. We would like to formulate this
as a vague conjecture.
\begin{conj}
{\rm 
For every variety $X$ with only Gorenstein toroidal singularities
there exists a loco sheaf $\msv(X)$ which is locally
isomorphic to the product of $\msv({\rm open~ball})$ 
and $\msv({\rm singularity})$ defined above. This construction
depends on the choice of parameters $g_n$ and perhaps 
on some other structures yet to be determined. 
Sheaf $\msv(X)$ is provided with the structure of a sheaf 
of conformal vertex algebras, and with $N=2$ structure
if $X$ is Calabi-Yau. The $L[0]=0$ component of this 
sheaf has a natural grading and differential which 
generalizes de Rham differential. The hypercohomology
of this complex has dimensions prescribed by
\cite{BD} and possesses a pure Hodge structure if $X$ is projective.
}
\label{msvgor}
\end{conj}

\section{Hypersurfaces in Gorenstein Toric Fano Varieties: general case}
Even though we are unable to construct \MSV\ for an arbitrary variety 
with Gorenstein toroidal singularities, the situation is somewhat better
in the special case of a hypersurface $X$ in a Gorenstein toric variety.
We can use the formulas for the smooth case applied now to
arbitrary cones. The resulting sheaf turns out to be loop-coherent. We
are mostly interested in the case when the ambient variety $\PP$ is 
Fano and the hypersurface $X$ is Calabi-Yau and generic, but most
statements hold true for any generic hypersurfaces. We will try
to extend the calculation of Sections 6-8.

 We  use the same notations $\Delta$, $\Delta^*$,
$\Sigma$, $M=M_1\oplus \ZZ$, $N=N_1\oplus \ZZ$. ${\rm deg}$, ${\rm
deg}^*$ as in Sections 7 and 8. We again consider a projective variety
$\PP$, a line bundle $L$ on it, and a hypersurface $X$ in $\PP$ given by
$f:\Delta\to\CC$. We define $\msv(X)$ as follows.

\begin{dfn}
{\rm
Let $f:\Delta\to \CC$ be a set of coefficients that defines $X$
and $g:\Delta^*\to \CC$ be a generic set of parameters. Then
for any cone $C^*\in\Sigma$ that contains $\deg^*$ sections of quasi-loco
sheaf $\msv(X)$
over the affine
chart $\AA_C$ are defined as BRST cohomology of $\Fock_{M\oplus C^*}$
with BRST operator
$${\cal BRST}_{f,g}=\oint(
\sum_{m\in \Delta} f_m
(m\cdot \Phi)(z) \ee^{\int
m\cdot B(z)}+
\sum_{n\in \Delta^*} g_n
(n\cdot \Psi)(z) \ee^{\int
n\cdot A(z)}
)dz.$$ 
}
\label{msvhyper}
\end{dfn}

For the above definition to make sense, we should show that the spaces
of sections constructed above are compatible with localization. Moreover,
we must show that they are loop-modules over the structure ring of $X$,
which means that they are annihilated by $f$.
\begin{prop} 
{\rm 
The above definition indeed defines a quasi-loco sheaf
of vertex algebras over $X$. It is provided with the structure of
topological algebra by formulas of Proposition \ref{jlgq}.
}
\label{indeed}
\end{prop}

{\em Proof.} To prove compatibility with localizations, we need to
show that for any cone $C^*$ of maximum dimension the $\BRST_{f,g}$
cohomology of $\Fock_{M\oplus C^*}$ is non-negatively graded with respect
to $L[0]$. Then the argument of Proposition \ref{localization} shows the
compatibility. The field $L[z]$ here is given by the formulas of
Proposition \ref{jlgq}, in particular, it differs slightly from $L(z)$ of
the Section 9. To avoid confusion we will call the operator given in
\ref{jlgq} by $L_{X,A}[0]$. This notation is chosen to indicate that we
are dealing with the Virasoro algebra of A model on the hypersurface $X$.
Explicitly,
$$
L_{X,A}(z)=\,\no B(z)\cdot A(z)\no +\no \partial_z\Phi(z)\cdot\Psi(z)\no
-{\rm deg}^*\cdot\partial_z A(z)
$$
so $L_{X,A}[0]$ counts the opposite of the sum of mode numbers of
$A,B,\Phi,\Psi$ plus $m\cdot n$ plus $m\cdot\deg^*$.
One can split $\BRST_{f,g}$ as a sum of $\BRST_f$ and $\BRST_g$ as usual.
Then we will have a spectral sequence as in Proposition \ref{affbfg}.
It is easily shown to be convergent because of Lemma \ref{lemma2}.
As a result, it is enough to show that the $\BRST_g$ cohomology 
of $\Fock_{M\oplus C^*}$ has nonnegative $L_{X,A}[0]$ eigenvalues.

Since $C^*$ contains $\deg^*$ and $\PP$ is Gorenstein, cone $C$ has some
special properties. One of the generators of its one-dimensional faces is
a vertex $m_0$ of $\Delta$, and all other generators lie in $M_1$. 
Notice that $\Fock_{M\oplus C^*}$ naturally splits as a tensor product of
$\Fock_{M_1\oplus C^*_1}$ and $\Fock_{\ZZ m_0\oplus \ZZ_{\geq 0}\deg^*}$.
Here $C_1^*$ is the cone in $N_1$ obtained by projecting $C^*$ there along
$\deg^*$. Moreover, it is easy to see that the $\BRST_g$ cohomology of 
$\Fock_{M\oplus C^*}$ is the tensor product of $\BRST_g$ cohomology of
$\Fock_{M_1\oplus C^*_1}$ and $\BRST_g$ cohomology of 
$\Fock_{\ZZ m_0\oplus \ZZ_{\geq 0}\deg^*}$. BRST operator on the first
space is defined precisely as in the Section 9, and BRST operator on the
second space is 
$$\BRST(z)=\deg^*\cdot \Psi (z)\ee^{\int \deg^*\cdot A(z)}.$$
Moreover, $L_{X,A}[0]$ is the sum of $L[0]$ from the Section 9 applied 
to $M_1\oplus C_1^*$ and 
$$
L_2[0]=\no(m_0\cdot B)(\deg^*\cdot A)\no[0] +
\no(\partial m_0\cdot \Phi)(\deg^*\cdot \Psi)\no[0]
+\deg^*\cdot A[0]
$$
Proposition \ref{boundedL} assures that BRST cohomology of $M_1\oplus
C_1^*$ doesn't have negative eigenvalues of $L_{X,A}[0]$. Explicit
calculation of BRST cohomology of $\Fock_{\ZZ m_0\oplus \ZZ_{\geq
0}\deg^*}$ given in Theorem \ref{dimone} allows us to conclude
that $L_2[0]$ also has non-negative eigenvalues. This assures that
$L_{X,A}[0]$ has only non-negative eigenvalues. 

However, we have only shown so far that $\msv(X)$ is a sheaf of vertex 
algebras over $\PP$. We also need to prove that the structure of
$\CC[C]$-loop-module induced from $\Fock_{M\oplus C^*}$ naturally gives
the structure of $\CC[C\cap M_1]/r_f$-loop-module, where $r_f$ is the
local equation of the hypersurface. It is enough to show this for a cone
$C^*$ of maximum dimension. Locally the element $r_f$ is
$$\sum_{m\in \Delta} f_m x^{m-m_0}.$$ 
Notice that the corresponding field 
$$r_f(z)=\sum_{m\in \Delta} f_m \ee^{(m-m_0)\cdot B(z)}$$
in $\Fock_{M\oplus C^*}$ could
be expressed as an anti-commutator of $\BRST_{f,g}$
and 
$$R(z)=\deg^*\cdot \Psi(z) \ee^{\int -m_0\cdot B(z)}.$$
Really, the OPE of this field with $\BRST_g(w)$ is non-singular, because
either $n\neq \deg^*$ and $-m_0\cdot n=0$, or $n=\deg^*$, which gives 
$-m_0\cdot n=-1$. However, in the latter case, we will also have $\deg^*
\Psi(z) \deg^* \Psi(w)$ which is $O(z-w)$, so overall the OPE is still
non-singular.  Hence, all modes of $r_f(z)$ act trivially on the 
cohomology by $\BRST_{f,g}$, which finishes the proof of the proposition.

We remark that it is plausible that all modes of $\ee^{m_0\cdot B(z)}$
act trivially as well, but we do not need to prove it, because 
$\CC[C\cap M_1]$ is embedded in $\CC[C]$.
$\hfill\Box$

It seems certain that $\msv(X)$ is loco with respect to $L_{X,A}[0]$,
but we do not have a proof of it yet. It would follow from any reasonable
solution of Conjecture \ref{msvgor}. On the other hand, we can easily  
show that
$\msv(X)$ is loco with respect to the grading $L_{X,B}[0]$ that comes from
the B model Virasoro field. This field is given by 
$$L_{X,B}(z)=\,\no B(z)\cdot A(z)\no -\no \Phi(z)\cdot
\partial\Psi(z)\no-{\rm deg}\cdot\partial_z B(z).
$$

\begin{theo}
{\rm
$\msv(X)$ is a loop-coherent sheaf with respect to the grading
$L_{X,B}[0]$.
}
\end{theo}

{\em Proof.}
The question is local and it is enough to consider a cone $C^*$
of maximum dimension. We are working in the set-up of the previous
proposition. By Theorem \ref{Bproploco}, cohomology of $\Fock_{M\oplus
C^*}$ with respect to $\BRST_g$ has graded components that are
Noetherian almost-modules over $\CC[C]$. Really, expressions for $L_B$ and
$L_{X,B}$ are identical (which is not the case for A model).

Notice that for a sufficiently big integer $k$, we can express 
$$\ee^{km_0\cdot B(z)}$$
as an anticommutator of some field and $\BRST_{f,g}$. That field is
similar to the one in the proof of Proposition \ref{homop}, but without
the extra $\ee^{-km_0\cdot B}$. The spectral sequence from
$\BRST_f$ cohomology of $\BRST_g$ cohomology to the
$\BRST_{f,g}$ cohomology degenerates by Lemma \ref{lemma2}, so
$\BRST_{f,g}$ cohomology is has graded components that are Noetherian
almost-modules over $\CC[C]/(x^{km_0})$. Since we have already shown
that $r_f$ acts trivially, these spaces are Noetherian over the structure
ring of $X$. $\hfill{\Box}$

\begin{coro}
{\rm For a fixed pair of eigen-values of $L[0]$ and $J[0]$, the
corresponding eigen-spaces of $H^*(\msv(X))$ are finite-dimensional.} 
\end{coro}

It is our firm belief that after Conjecture \ref{msvgor} is
successfully proved, the sheaf $\msv(X)$ could be identified as a
(generalized) chiral de Rham complex of $X$ as implied by this notation.
However, it is still well defined as a sheaf of vertex algebras and one
may ask how to calculate its cohomology. 

We can also ask whether the analog of Proposition \ref{affbfg} 
still holds. For orbifold singularities Proposition \ref{affbfg} still
holds the way it is stated, but the proof must be different, because
vertical cohomology of the double complex considered there is non-zero for
more than one row. In the case of orbifold singularities spectral
sequence still degenerates, because we may  split the picture according to
the eigen-values of $B[0]$ modulo the lattice spanned by generators of
one-dimensional faces of $C^*$. This spectral sequence might degenerate
for every $C^*$, but we can not prove it. 

Unfortunately, Proposition \ref{zerototal} does not hold even for 
orbifold singularities. It is also not clear that double complex
of Theorem \ref{conebig} gives degenerate spectral sequences for 
arbitrary Gorenstein toric Fano varieties. On the other hand,
$\BRST_{f,g}$ cohomology of $\Fock_{M\oplus K^*}$ could still be
the correct vertex algebra to consider, once the relation to physicists'
A and B models becomes more clear. Section 8 never uses the fact that
$\PP$ is non-singular and generalizes to any $\PP$. To sum it up, it is
plausible that Theorem \ref{maintheo} holds for any toric Gorenstein Fano
varieties, but we can only prove it in the smooth case.

To complete the discussion we must mention that Mirror Symmetry
for Calabi-Yau complete intersections in Gorenstein toric
Fano varieties can be adequately treated by the methods 
of this paper. It is appropriate to state the final
conjecture that covers all examples of "toric" Mirror Symmetry.
See \cite{BBci} for relation between complete intersection
examples of Mirror Symmetry and pairs of dual reflexive Gorenstein
cones.

\begin{conj}
{\rm
Let $M$ and $N$ be dual lattices with dual cones $K$ and $K^*$
in them. We assume that $K$ and $K^*$ are reflexive Gorenstein,
which means that $K\oplus K^*$ is Gorenstein in $M\oplus N$.
We denote the corresponding degree elements by ${\rm deg}$
and ${\rm deg}^*$. We are also provided with generic numbers
$g_n$ and $f_m$ for all elements in $K$ and $K^*$.
Then, if reflexive cones come from Calabi-Yau complete intersections,
 vertex algebras of these Calabi-Yau manifolds are degenerations
of BRST cohomology  of 
$\Fock_{M\oplus N}$ by operator
$$\Bfg=\oint(
\sum_m f_m(m\cdot \Phi)(z) \ee^{\int m\cdot B(z)}+
\sum_n g_n(n\cdot \Psi)(z) \ee^{\int n\cdot
A(z)})dz.$$
The degeneration is provided by fans that define the corresponding
toric varieties.
Structure of topological algebras of dimension $\dim M-2deg\cdot
deg^*$ is given by
$$
Q(z)=A(z)\cdot\Phi(z)-{\rm deg}\cdot\partial_z \Phi(z),~
G(z)= B(z)\cdot\Psi(z)-{\rm deg}^*\cdot\partial_z \Psi(z),
$$
$$
J(z)=\,\no\Phi(z)\cdot\Psi(z)\no + \,{\rm deg}\cdot B(z)
-{\rm deg}^*\cdot A(z),  
$$
$$
L(z)=\,\no B(z)\cdot A(z)\no +\no \partial_z\Phi(z)\cdot\Psi(z)\no
-\,{\rm deg}^*\cdot\partial_z A(z).
$$
}
\label{ci}
\end{conj}

\section{Open questions and concluding remarks}
In this section we point out important questions that 
were not addressed in this paper as well as possible applications
of our results and techniques.

\noindent $\bullet$ 
%Flatness
It remains to show that deformations of the Master Family
of vertex algebras are flat in the appropriate sense.
For instance, we would love to say that dimensions of 
$L[0]$ eigen-spaces are preserved under these deformations.

\noindent $\bullet$ 
%new conformal theories
One can generalize the construction of Conjecture
\ref{ci} to go from $K\oplus K^*$ in $M\oplus N$ to
any Gorenstein self-dual cone in a lattice with
inner product by using the BRST field
$${\cal BRST}(z)= {\sum_n}g_n (n\cdot {\rm fermion})(z)
\ee^{\int n\cdot {\rm boson}(Z)}.$$
These theories still have conformal structure,
with $L$ given as $L_{\rm flat} - {\rm deg}\cdot\partial {\rm boson}$.
Do these theories have any nice properties or physical significance?

\noindent $\bullet$ 
%Unitarity
It would be great to provide vertex algebras of Mirror Symmetry
with unitary structure. Inequalities on eigen-values of
$L[0]$ and $J[0]$ seem to suggest its existence.

\noindent $\bullet$
%Correlators
It is extremely important to use results of this paper to
get actual correlators of corresponding conformal field theories
and thus to draw a connection to the calculation of \cite{candelasetal}.

\noindent $\bullet$ 
%Kontsevich
There must be a connection between calculations of this paper
and moduli spaces of stable maps defined by Kontsevich. It
remains a mystery at this time.

\noindent $\bullet$
%GKZ
It would be interesting to see how GKZ hypergeometric system
enters into our picture. Solutions of GKZ system are known
to give cohomology of Calabi-Yau hypersurfaces, see for
example \cite{hosono,stienstra}.

\noindent $\bullet$ 
%Antiholomorphic part
One must also do something about the antiholomorphic part of the
$N=(2,2)$ super-conformal algebra. Perhaps this problem is not 
too hard and has its roots in the author's ignorance.

\noindent $\bullet$ 
%String-theoretic cohomology
One should construct generalized \MSV es as suggested in
Conjecture \ref{msvgor}. This seems to be a realistic project 
since all we really need to do is to extend the automorphisms
of toroidal singularities to the suggested local descriptions
of \MSV es. It is also interesting to see which of the 
standard properties of cohomology of smooth varieties
generalize to string-theoretic cohomology. 
 
\noindent $\bullet$ 
%Other Gorenstein singularities
One should define string cohomology for all Gorenstein, and
perhaps $\QQ$-Gorenstein singularities. See \cite{Batyrevnew}
for the definition of string cohomology numbers in this generality.
This may also shed some light on generalized McKay correspondence.

\noindent $\bullet$ 
%Possible application to hyperbolicity
Our results and techniques may have applications to hyperbolicity.
Indeed, for a smooth $X$ there is a loco subsheaf 
$\msv_{b,\varphi}$
of $\msv(X)$ which
is generated by modes of $b^i$ and $\varphi^i$ only. This part
of $\msv(X)$ is contravariant, so for any map from a line 
to $X$ it could be pulled back to it. Then global sections
of $\msv_{b,\varphi}$ could give restrictions on possible maps 
to $X$. On the other hand, rich structure of the whole 
$\msv(X)$ might help to show that there are plenty of such sections.

\noindent $\bullet$ 
%modular forms 
Finally, cohomology of $\msv(X)$ is graded
by $L[0]$ and $J[0]$, and one can show that ${\rm
Trace}(q^{L[0]}w^{J[0]})$ has some modular properties. It is 
directly related to elliptic genus of $X$. This issue will
be addressed in the upcoming joint paper with Anatoly Libgober.

\end{document}